

\input amstex \documentstyle{amsppt}

\NoRunningHeads
\magnification=1100
\baselineskip=29pt
\parskip 9pt
\pagewidth{5.2in}
\pageheight{7.2in}

\TagsOnLeft
\CenteredTagsOnSplits

\let\al=\alpha

\def\aln{\al_1\cdotss\al_n}
\def\alndh{\al_1\cdotss\al_n}
\def\aotilalg{A_1X/\!\sim_{\text{alg}}}

\let\be=\beta
\def\bl{\bigl(}
\def\br{\bigr)}
\def\Bl{\Bigl(}
\def\Br{\Bigr)}
\def\bem{\be_1\cdotss\be_m}

\def\cdotss{,\cdots,}

\def\c{\Cal{C}}
\def\cbnm{C_B^{n+m}}
\def\cbnmoo{C_{B_1}^{1+n_1+m_1}}
\def\cbnmto{C_{B_2}^{1+n_2+m_2}}
\def\cbnmo{C_{B_1}^{n_1+m_1}}

\def\codim{\operatorname{codim}}
\def\cb{\Cal C_B}
\def\cbo{C_{B_1}}

\def\Del{\Delta}
\def\dual{^{\vee}}
\def\del{\delta}

\def\doublepri{^{\prime\prime}}

\let\eps=\epsilon
\let\ep=\epsilon
\def\epmdh{\ep_1\cdotss\ep_m}
\let\epm=\epmdh

\def\Ext{\operatorname{Ext}}

\def\f{\Cal F}

\def\char{\text{char}\,}

\def\hugeindex{_{\Sb B_1+B_2=B,\, m_1+m_2=m,\\\tau\in S(m_1,m)\endSb}}
\def\h{\Cal H}

\def\Hom{\operatorname{Hom}}

\def\indexn{_{i=1}^n}
\def\indexm{_{j=1}^m}

\def\Int{\text{Int}_B}
\def\int{\text{Int}}
\def\Im{\operatorname{Im}}

\def\lsta{_{\ast}}
\def\Lam{\Lambda}

\def\lra{\longrightarrow}
\def\length{\text{length}}

\def\morsb{\text{Mor}(\sig,B)}
\def\mh{\!:\!}
\def\mapright#1{\,\smash{\mathop{\lra}\limits^{#1}}\,}

\def\Mor{\text{Mor}}
\let\midd=\mid

\def\OO{\Cal O}
\def\opcap{\operatornamewithlimits{\bigcap}}
\def\opprod{\prod}
\def\opcup{\operatornamewithlimits{\bigcup}}

\def\pri{^\prime}
\def\piis{\pi_{i\sig}}
\def\piix{\pi_{iX}}
\def\piso{p_{i\sigo}}
\def\pix{p_{iX}}

\def\qix{q_{iX}}
\def\qist{q_{i\sigt}}

\def\px{p_X}
\def\ps{p_{\sig}}

\def\phinm{\varPhi_B}
\def\psinm{\Psi_B}

\def\piiz{\pi_{iZ}}
\def\Pnm{P^{n|m}}

\def\r{\Cal R}

\def\siui{S_i^{\upsilon_i}}

\def\stimex{\sig\time X}
\def\stimev{\sig\time V}
\def\s{\Cal S}
\def\spec{\operatorname{Spec}}
\def\skl{s_1\cdotss s_{n+m}}
\def\sixnotm{\bl\sig\time X\br^{n_1}\time \bl \sig\time X)^{n_2}\time X^m}
\let\si=\sigma

\def\sindh{\si_1\cdotss\si_{n}}

\def\sn{s_1\cdotss s_n}

\let\sig=\Sigma
\def\sigo{\sig_1}
\def\sigt{\sig_2}
\def\sta{^{\ast}}
\def\sigox{\sig_1\time X}

\def\stimex{\sig\times X}
\def\sub{\subseteq}
\def\sta{^{\ast}}
\def\sixnm{\bl \sig\time X\br^n\time X^m}

\def\tilzeta{\tilde\zeta}
\def\tilF{\widetilde{F}}

\def\tIntnotm{\widetilde{\text{Int}}^{n_1:n_2}_B}
\def\tInt{\widetilde{\text{Int}}}
\def\time{\!\times\!}
\def\topo{^{\circ}}
\def\tcbnm{\c_B^{n+m}}
\def\tilzeta{\tilde\zeta}
\def\tkls{\Gamma(k_1,\Lam_1,\upsilon)}
\def\t{\Cal T}

\def\varphis{\varphi_{\sig}}
\def\varphix{\varphi_X}
\let\vphi=\varPhi
\let\vpsi=\varPsi

\let\ups=\upsilon

\def\umo{^{-1}}
\def\udot{^{\,\cdot}}
\def\ucirc{^{\circ}}
\def\usi{^{\si}}

\def\wb{\Cal W_B}
\def\w{\Cal W}
\def\wbcirc{W_B^{\circ}}

\def\xn{x_1\cdotss x_n}
\def\xx{\chi}

\def\Yn{Y_1\cdotss Y_n}

\def\yisi{Y_i^{\si_i}}
\def\Ydot{Y_{\cdot}}

\def\z{\Cal Z}
\def\zndh{z_1\cdotss z_n}
\def\zn{z_1\cdotss z_n}
\def\Zm{Z_1\cdotss Z_m}
\def\zdot{z_{\cdot}}
\def\ZZ{{\Bbb Z}}

\def\Zdot{Z_{\cdot}}
\def\zjej{Z_j^{\eps_j}}

\let\pro=\proclaim
\let\endpro=\endproclaim
\def\proof{{\it Proof.\ }}


\topmatter
\title
The quantum cohomology of homogeneous varieties
\endtitle

\author
Jun Li and Gang Tian
\endauthor
\thanks
The first named author was partially supported
by NSF grant DMS-9307892 and A. Sloan fellowship
and the second named author was
partially supported by NSF grant DMS-9303999 and Alan T.
Waterman award.
\endthanks

\affil
Mathematics Department\\
University of California, Los Angeles\\
and\\
Courant Institute of Mathematics\\
New York University
\endaffil

\address
Mathematics Department,
University of California, Los Angeles, CA 90024
\endaddress

\email
jli@math.ucla.edu
\endemail

\endtopmatter
\document

\head 0. Introduction
\endhead

The notion of quantum cohomology was first proposed by Witten
[Va, Wi], based on topological field theory. Its mathematical theory
was only established recently by Y. Ruan and the second named author
[RT, Ru], where they proved the existence of the quantum rings on
semi-positive symplectic manifolds. (Fano manifolds are particular
semi-positive manifolds.) In this note, we will provide a purely
algebro-geometric proof to the existence of quantum cohomology
rings for a special class of manifolds already treated there, namely
homogeneous manifolds. By using algebro-geometric approach, we can
prove the existence of quantum cohomology of homogeneous
varieties defined over any algebraically closed field.
This should be useful in enumerative geometry. We believe the
approach here can be applied to a larger class of algebraic varieties,
such as toric varieties.

Let $X$\ be a $d$-dimensional smooth variety. To construct its quantum
cohomology of $X$, one needs to define a class of
enumerative invariants, which
is the Gromov-Witten invariants when $X$\ is a complex
variety. Let $\sig_g$\
be the smooth curve of genus $g$\ with marked points $\xn\in\sig_g$, let
$B\in \aotilalg$\ and let $\Mor(\sig_g,B)$\ be the moduli
scheme of morphisms $f\mh \sig_g\to X$\ satisfying $[f(\sig_g)]\in B$.
For any cohomology classes $\aln\in A\sta X$\ whose Poincare duals
are represented by subvarieties $Y_1\cdotss Y_j\sub X$,
the enumerative invariant $\vphi_{(B,g)}$\ is defined as
$$\vphi_{(B,g)}(\aln)=\sum_{f\in\Lam}m(\Lam,f)
\tag 0.1
$$
when $\dim \Lam=0$\ and is zero otherwise,
where $\Lam$\ is the scheme
$$\Lam=\{f\in\Mor(\sig_g,B)\mid f(x_i)\in Y_i\}
$$
and $m(\Lam,f)$\ is the multiplicity of $\Lam$\ at $f$. Clearly,
$\vphi_{(B,g)}$\ is symmetric. By Riemann-Roch
theorem, (0.1) should be non-trivial only when
$$\sum_{i=1}^n\codim Y_i-c_1(X)\cdot B+d(1-g)=0,
\tag 0.2
$$
since its left hand side is the expected dimension of $\Lam$.

For the moment, let us assume the enumerative invariant
$\vphi_{(B,g)}$\ has been defined
for all possible $(g, B,\alndh)$. Then they should obey a simple
composition law [RT]:
Assume that the dual of the diagonal cycle $\Del\sub X\times X$\ has
the following Kunneth decomposition
$$[\Del]\dual=\sum_{l=1}^p\zeta_l\times\tilzeta_l
$$
(via obvious inclusion $A\sta X\times A\sta X\sub A\sta(X\times X)$),
then for any partition $g=g_1+g_2$\ with $g_1,g_2\geq 0$\ and
$n_1$\ lies between $\max(0, 2-2g_1)$\ and $\min(n,n+2g_2-2)$,
$$\vphi_{(B,g)}(\aln)=\sum_{\Sb B=B_1+B_2\\ 1\leq l\leq p\endSb}
 \vphi_{(B_1,g_1)}(\zeta_l,\al_1\cdotss\al_{n_1})\cdot
 \vphi_{(B_2,g_2)}(\tilzeta_l,\al_{n_1+1}\cdotss\al_n).
$$
The invariants $\vphi_{(B,0)}$\ are used to define a new ring structure
on $A\sta X$, called the quantum ring, whose importance in
enumerative geometry is yet to be realized. This in part is due to the
lack of rigorous algebraic definition, except the case treated in [CM,
Be] and this note.

This new ring structure on $A\sta X$\ is defined as follows:
We first form a formal infinite series
$$\widetilde{\vphi}(\al,\be,\gamma)=\sum_{B\in\aotilalg}
\vphi_{(B,0)}(\al,\be,\gamma)\, e^{-t(c_1(X)\cdot B)},
$$
where $\al,\be,\gamma\in A\sta X$\ and $(c_1(X)\cdot B)$\ is the
the degree of $c_1(X)(B)\in A_0X$. Then we define a new product
$$\times_{\bold Q}: A\sta X\times A\sta X\lra A\sta X
$$
by the rule
$$\al\times_{\bold Q}\be=\sum_{l=1}^p \widetilde{\vphi}(\al,\be,
\tilzeta_l)\,\zeta_l.
$$
This product gives rise to a ring structure (the only non-trivial
part is to check its associativity) that follows immediately from the
composition law above. It is worth mentioning that some enumerative
problems, say counting curves in Fano varieties, can be worked
out by using the homomorphism $\vphi_{(B,g)}$\ once
the composition law is established. As Witten pointed out,
a simple corollary of composition law is the associativity law for
rational curves.
For instance, assuming the associativity law,
Kontsevich calculated the number of degree $d$ rational
curves in $CP^2$ through $3d-1$ generic points.
One can also deduce an explicit formula for
counting rational curves in $\bold P^n$\ by using the
associativity law for rational curves
[KM, RT].

There are several issues that need
to be addressed before this set of invariants can be defined.
First, $\Lam$\ has to be discrete under (0.2) in order to make sense of
(0.1), at least not appealing to the excessive intersection
theory. Usually, this rarely happens
even after moving $Y_i$'s to general positions (there are
limitations in ``moving'' subvarieties). In [RT], this
difficulty was overcome by replacing $f$\ with solutions of the
Cauchy-Riemann equation with respect to a generic almost
complex structure on $X$. In this note, we
choose to avoid this difficulty by restricting ourselves to the case
where $g=0$\ and $X$\ is a
homogeneous variety (strictly speaking we will only consider those
homogeneous varieties whose tangent bundles are generated by
global sections), since then the space of morphisms $\text{Mor}(
\sig_0,B)$\ is smooth everywhere and the variety $Y_i$\
can be made in a general position after translation via a
$\si\in G$\ [Kl]. Another issue is that the space
$\text{Mor}(\sig_0,B)$\ is never proper. This can be dealt
easily with for homogeneous variety by using non-proper intersection
[Fu, \S10.2]. The main input of this note is the construction of a degeneration
of $\text{Mor}(\sig_0,B)$\ upon which the composition law follows
immediately.

Before sketching the construction of this degeneration, let us first
state the main theorem of this paper. In the following, for any variety
$X$\ we denote by ($A^iX$) $A_iX$\ the birational (co)homology group
of $X$\ [Fu, \S17.3]. In case $\si$\ is
a transformation of $X$, then we denote by $Y\usi$\ the translation
$\si(Y)\sub X$. We denote by $\sig$\ the smooth rational curve.

\pro{Definition-Theorem 0.1}
Let $K=\overline K$\ be an algebraic closed filed and let $X=G/P$\ be
a $d$-dimensional smooth projective variety over $K$\ that is a
quotient of a linear algebraic group $G$\ by a closed subgroup $P$.
Then for any $B\in \aotilalg$, there is a (group) homomorphism
$$\phinm(\cdot\mid\cdot):
\prod\indexn A\sta X\times\prod\indexm A\sta X\lra \ZZ
$$
of which the following holds: Let $\al_i\in A^{d_i}X$, $i=1\cdotss n$,
and $\be_j\in A^{e_j}X$, $j=1\cdotss m$, be any cohomology classes.
Then $\phinm\bl\aln\midd\bem\br=0$\ unless
$$\sum\indexn d_i+\sum\indexm e_i=c_1(X)\cdot B+d+m.
$$
In this case if we further assume the Poincare dual of
$\al_i$\ (resp. $\be_j$) are represented by subvariety $Y_i\sub X$\
(resp. $Z_j\sub X$), then for general
$\zndh\in\sig$\ and general $\sindh,\epmdh\in G$,
$$\phinm\bl\aln\midd\bem\br=\sum_{f\in\Lam}m(\Lam,f),
$$
where $\Lam$\ is the scheme
$$\Lam=\left\{(f,w_1\cdotss w_m)\Biggm\vert
\matrix f\in\morsb,\ w_1\cdotss w_m\in\sig\\
f(z_i)\in\yisi,\ 1\leq i\leq n;\
f(w_j)\in Z_j^{\ep_j},\  1\leq j\leq  m.
\endmatrix
\right\}
\tag 0.3
$$
(that is discrete) and $m(\Lam,f)$\ is the multiplicity of $\Lam$\ at $f$.
\endpro

Here by a property holds for general $\zn\in\sig$\ and
$\sindh\in G$\ we mean that there is a dense
open sebset $U\sub \sig^n\times G^n$\ such that the property holds for all
$(\zn,\sindh)$\ in $U$.

\pro{Theorem 0.2}
Let the notation be as in definition-theorem 0.1. Then the enumerative
invariant $\phinm$\ satisfies the following composition law:
Let $n=n_1+n_2$, $n_1,n_2\geq 2$, be any partition and let $[\Del]\dual
=\sum_{l=1}^p\zeta_l\times\tilzeta_l$\ be the Kunneth decomposition
\footnote{It is shown in [FMS$^2$] that such decompositions always exist for
the class of varieties considered in this paper.}
of the Poincare dual of the diagonal $\Del\sub X\times X$, then
$$\align
&\phinm(\aln\mid\bem)\\
=&\sum_{\Sb B=B_1+B_2,\,\tau\in S(m)\\ 1\leq l\leq p,\,0\leq k\leq m
\endSb}
{1\over k!(m-k)!} \cdot
  \vphi_{B_1}(\zeta_l,\al_1\cdotss\al_{n_1}\mid\be_{\tau(1)}
\cdotss\al_{\tau(k)})\\
&\qquad\qquad\qquad\qquad\qquad\cdot
\vphi_{B_2}(\tilzeta_l,\al_{n_1+1}\cdotss\al_{n}\mid
\be_{\tau(k+1)}\cdotss\be_{\tau(m)}).
\endalign
$$
Here $S(m)$\ is the group of permutations of $\{1,\cdots,m\}$.
\endpro

For any smooth Fano variety Y. Ruan and the second named author
defined the mixed quantum invariant
$$\vphi_{B\pri}^{\text{Sym}}: \prod_{i=1}^nH_{\ast}(X,\ZZ)\times
\prod_{j=1}^m H_{\ast}(X,\ZZ)\lra \ZZ,
$$
where $B\pri\in H_2(X)$,
by analytic method. The homomorphism defined in this paper is
compatible to $\vphi_{B\pri}^{\text{Sym}}$\
for the class of varieties considered in this paper.

\pro{Theorem 0.3}
Let $X$\ be any smooth complex projective variety acted on
transitively by a linear algebraic group and let
$$\mu: A^p X \lra H_{2p}(X,\ZZ)
$$
be the homomorphism that is the composition of the Poincare
dual $A^p X\to A_pX$\ and the cycle map $cl\mh A_pX\to H_{2p}(X)$\
[Fu, \S9.1]. Then for any $\aln,\bem\in A\sta X$,
$$\vphi_B\bl\aln\mid\bem\br=
\vphi_{B\pri}^{\text{Sym}}\bl\mu(\al_1)\cdotss\mu(\al_n)\mid
\mu(\be_1)\cdotss\mu(\be_m)\br,
$$
where $B\pri=cl(B)$.
\endpro

\pro{Remark}
In general, even when $\vphi_B$\ can be defined as the degree of
a 0-cycle (as in (0.3)), it is not obvious that it will coincide with
$\vphi_{B\pri}^{\text{Sym}}$, originally defined as a symplectic
invariant.
\endpro

Now we sketch the proof of Theorem 0.2. For simplicity, we will
consider the case $m=0$. After identifying $f\in\morsb$\ with its graph
$\Gamma_f\sub \sig\times X$, we can embed $\morsb$\ as a dense open
subset in $W_B$\ that consists of subschemes $C\sub\sig\time X$\
satisfying
$$\dim C=1, \xx(\OO_C)=1, [p_X(C)]\in B\ \text{and}\,
\deg(C,\ps\sta\OO_{\sig}(1))=1,
$$
where $\ps$\ and $\px$\ are projections from $\sig\times X$\ to
its factors. (Strictly speaking, $W_B$\ is the closure
of $\morsb$\ in the respective Hilbert scheme.)
$W_B$\ is smooth at $\Mor(\sig,B)$\ since $TX$\ is generated by
global sections by our assumption on $P\sub G$.
Let $C_B$\ be the universal family and let $C_B^n$\ be the Cartesian product
of $n$\ copies of $C_B$\ over $W_B$.
For any $z\in\sig$\
and subvariety $Y\sub X$, we can define an incidence subscheme
$$\Lam_i(z,Y)=\left\{(C,x_1\cdotss x_n)\Bigm\vert
\matrix C\in W_B,\ x_1\cdotss x_n\in C,\\
\ps(x_i)=z\ \text{and}\ \px(x_i)\in Y_i
\endmatrix
\right\}.
\tag 0.4
$$
By using the general position theorem
for homogeneous varieties (see Lemma 2.3),
for any $Y_1,\cdots,Y_n\sub X$\ satisfying
$\sum\codim Y_i=c_1(X)\cdot B+d$,
we can find general $z_1,\cdots,z_n\in \sig$\ and general
$\si_1\cdotss\si_n\in G$\ such that the translations
$\yisi$\ of $Y_i$\ are in general position in the sense that
$$\opcap_{i=1}^n \Lam_i(z_i,\yisi)\sub\morsb\sub W_B
\tag 0.5
$$
and is discrete and Cohen-Macaulay. Thus $\vphi_B([Y_1]\dual\cdotss
[Y_n]\dual)$\ can be defined
as the degree of the zero cycle (0.5).

To establish the composition law, we need to construct a degeneration of
$W_B$. We take a simple minded approach: Let $0\in V$\ be a smooth curve
as a parameter space and let $Z$\ be the blowing-up of $\sig\times V$\
at a point over $0\in V$. $Z$\ is a degeneration of $\sig$\ into a nodal
curve that has two components, say $\sigo$\ and $\sigt$.
We let $\wb$\ be the relative Hilbert scheme of
one dimensional subschemes $C\sub Z\times X$\ obeying constrains
similar to that of $W_B$\ and let $\cb$\ and $\cb^n$\ be the universal
family and the product of this family over $\wb$\ as $C_B$\ and $C_B^n$\ do
(see \S 3 for details).
Again, for any section $z$\ of $Z\to V$\ and subvariety $Y\sub X$, we can
define an incidence subscheme $\tilde\Lam_i(z,Y)\sub\cb^n$\ as
$$\tilde\Lam_i(z,Y)=\left\{(S,\sn)\in\cb^n\Bigm\vert
\matrix S\in\wb,\ \sn\in S,\\
p_Z(s_i)\in z\ \text{and}\ \px(s_i)\in Y
\endmatrix
\right\}.
$$
Now let $z_{1}\cdotss z_{n}$\ be $n$\ sections of $Z\to V$\ of which the first
$n_1$\ sections passing through the component $\sigo\sub Z_0$\ and the
remainder
passing through $\sigt\sub Z_0$. Since both $\sig$\ and $X$\ are
homogeneous, for the $Y_1\cdotss Y_n\sub X$\ as before, we can choose
$z_{i}$\ and $\si_i\in G$\ general so that the intersection
$$\opcap_{i=1}^n \tilde\Lam_i\bl z_{i},\yisi\br\sub \cb^n
\tag 0.6
$$
is flat, finite and Cohen-Macaulay over $V$\ at $0\in V$. Hence
$$\vphi_B([Y_1]\dual\cdotss[Y_n]\dual)=\deg\Bigl(\bl
\opcap_{i=1}^n\tilde\Lam_i(z_{i},\yisi)\br\bigcap\c^n_{B,t}\Bigr)
\tag 0.7
$$
for all $t\in V$\ near $0$,
where $\c^n_{B,t}$\ is the fiber of $\c^n_B$\ over $t\in V$.
Now it is straight forward to check that the points of the intersection
in (0.6) over $0\in V$\ are pairs of morphisms
$f_1,f_2\mh \sig\to X$\ satisfying properties resemble that of (0.5).
Thus (0.7) can be counted based on intersection of
incidence schemes in $C_{B_1}^{\,\cdot}\times C_{B_2}^{\,\cdot}$, where
$B_1$\ and $B_2$\ runs through all
possible $B_1+B_2=B$. This line of argument eventually
will lead to the proof of the composition law as
stated in theorem 0.2.

We should say one word about the assumption on the closedness of
$P\sub G$. When $\char K=0$, it is equivalent to assuming $G$\ is a
linear algebraic group acted on transitively on a smooth projective
$X$. In case $\char K>0$, this condition can be relaxed by in addition
to the previous assumption on $G$\ and $X$, we require that for any
morphism $f\mh \sig\to X$\ $f\sta T_X\otimes \OO_{\Sigma}(1)$\ is ample. This
will ensure the smoothness of $\morsb$\ that is all we need in
proving the composition law.

While this paper was in writing, A. Bertram [Be] and B. Crauder and
R. Miranda [CM] have finished their works attacking the same problem
for some classes of varieties by algebraic approach. In [Be], A. Bertram
established the composition law for all Grassmannians by generalizing the
classical Giambelli and Pieri's formulas.
B. Crauder and R. Miranda in [CM] studied the
quantum rings of rational surfaces in detail, establishing the
composition law in the mean time. (I. Ciocan-Fontanine
has informed us that he has constructed the quantum cohomology rings of
flag varieties using algebraic geometry.)
We feel that the current approach is more direct, and hopefully
can be modified to prove the composition law for all smooth
Fano varieties, after their enumerative invariants have been defined
appropriately.

Most part of this work was completed while the second author was
visiting Stanford University in early 1994, for which he is grateful to
the host for prviding a wonderful environment.
We also like to thank W. Fulton for
helpful conversation and advice.

The layout of this paper is as follows: In \S 1, we will give the definition
of the the enumerative invariants for projective varieties.
\S 2 concerns rational curves in
homogeneous variety. Finally, we will prove the composition law
in \S 3.

\head 1. The definition of the quantum invariants
\endhead

In this section, we will first define the homomorphism $\vphi_B$\
by pure cohomological calculation on some auxiliary spaces.
Then we will show that in the ideal situation, its value
is exactly the
intersection of some ``geometrically defined'' cycles in the space of
morphisms from the rational curve $\sig$\ to $X$, as proposed by
physicists. In the end, we will compare our definition with the
Symplectic invariant defined in [RT].

We begin with some words on the convention that will be used
throughout this paper. Let $K=\overline K$\ be an algebraically
closed field. Let $X=G/P$\ be a $d$-dimensional
smooth projective variety and let $\sig$\ be the
smooth rational curve, both defined over $K$.
For any scheme $Y$, we denote by ($A^pY$)
$A_pY$\ the p-th (co)homology group of $Y$. By definition, any
element $\al\in A^pY$\ is a collection of homomorphisms
$A_kY'\to A_{k-p}Y'$, for all $Y'\to Y$\ and all $k$, satisfying certain
properties (see [Fu]).
Let $\ps$\ and $\px$\ be projections from $\sig\times X$\ to $\sig$ and
$X$\ respectively.

We first embed the scheme $\morsb$\ into some projective scheme,
since the former usually is not complete. The way
we treat this problem is to embed it into the
Hilbert scheme of subschemes of $\sig\times X$\ by identifying
a morphism with its graph.
For any $B\in\aotilalg$, we form the Hilbert scheme
$H_B$\ of subschemes $C\sub \stimex$\ satisfying the following
constrains:
$$\dim C=1,\ \xx(\OO_C)=1,\ [\ps(C)]\in B\ \text{and}\
\deg(C,\ps\sta(\OO_{\sig}(1))=1.
\tag 1.1
$$
By definition, each point of $H_B$\ corresponds to a complete
1-dimensional subscheme $C\sub\stimex$, which in the
sequel will be denoted by
$C\in H_B$. Let $\varphis\mh C\to \sig$\  and
$\varphix\mh C\to X$\
be morphisms induced by projections $\ps$\ and $\px$.
When $C$\ is irreducible, $\deg(C,\varphis\sta(\OO_{\sig}(1))=1$\
implies that at general points, $C$\ is reduced and
$\varphis\mh C\to \sig$\ is one-one. Then because
$\xx(\OO_C)=1$, $\varphis$\ is an isomorphism.
Hence $C$\ gives rise to a morphism
$$f_C=\varphix\circ\varphis\umo: \sig\lra X
$$
with $[f_C(\sig)]\in B$. Also, by Riemann-Roch theorem, the expected
dimension of $H_B$\ at $C$\ is
$$\rho(B)=c_1(X)\cdot B+d.
\tag 1.2
$$
Now let $\morsb$\ be the moduli scheme of morphisms $f\mh \sig\to X$\
with $f\lsta([\sig])\in B$. By identifying $f$\ with its graph $\Gamma_f\sub
\stimex$, we obtain an open immersion $\morsb\sub H_B$. We let
$W_B^{\circ}$\ be the image of $\morsb$\ and let $W_B$\
be the scheme theoretic closure of $\wbcirc\sub H_B$.
Note that $\wbcirc\sub H_B$\ is the open subset of points
whose corresponding subschemes $C\sub\sig\time X$\ are irreducible.

For technical reason, in the remainder of this section we will assume
$W_B$\ has pure dimension $\rho(B)$, which is the case when $X$\ is
homogeneous. Let $C_B\sub\stimex\times W_B$\
be the universal family, flat over $W_B$, and let $\cbnm$\ be the Cartesian
product
$$\cbnm= \underbrace{C_B\times_{W_B}\cdots\times_{W_B}C_B}_
{n+m\ \text{copies}},
$$
where $n, m\geq 0$. Note that $\dim\cbnm=\rho(B)+n+m$\ and points
of $\cbnm$\ are tuples $(C,x_1\cdotss x_{n+m})$\ with
$C\in W_B$\ and $x_i\in C$.
We let $\pi_i\mh\cbnm\to C_B$\ be the projection
onto the i-th copy and let $\piis=\ps\circ\pi_i$\ and $\piix=\px\circ\pi_i$\
be compositions. We next define a projection
$$P^{n|m}:\cbnm\lra \bl\sig\times X\br^n\times X^m
\tag 1.3
$$
that is the product
$$\Bl\prod\indexn\piis\times\piix\Br\times
\Bl\prod\indexm\pi_{n+jX}\Br.
$$
By letting $e_0\in A^1\sig$\ be the Poincare dual of a point in
$\sig$, we get an inclusion
$$\prod\indexn A\sta X \time\prod\indexm A\sta X
\lra A\sta\bl(\sig\time X)^n\time X^m\br
\tag 1.4
$$
that sends
$$(\aln;\bem)\in
\prod\indexn A\sta X \times\prod\indexm A\sta X
$$
to $(\prod^ne_0\time\al_i)\time(\prod^m\be_j)$. We then define
$$\psinm(\cdot\midd\cdot): \opprod\indexn A\sta X\times\prod\indexm A\sta X
\lra A\sta\cbnm
\tag 1.5
$$
be the homomorphism that is the composition of the inclusion (1.4) with
the pull-back $\bl P^{n|m}\br\sta$. Alternatively, $\psinm$\ can be
defined as
$$\psinm(\aln\midd\bem)=
\Bigl(\opcup\indexn\piis\sta e_0\cup\piix\sta\al_i\Bigr)
\bigcup\Bigl(\opcup\indexm\pi_{n+jX}\sta\be_i\Bigr),
$$
where $\cup$\ is the cup-product.
By the definition of cohomology theory,
$$\psinm(\aln\midd\bem)([\cbnm])
\tag 1.6
$$
is an element in $A\lsta\cbnm$, where
$[C_B^{n+m}]\in A\lsta\cbnm$\ is the fundamental class. We then define
$$\phinm(\cdot\midd\cdot):\opprod\indexn A\sta X
\times\prod\indexm A\sta X \lra \ZZ
\tag 1.7
$$
to be the degree of the component of (1.6) in $A_0\cbnm$. Namely,
$$\phinm(\aln\midd\bem)=\deg\Bigl(\psinm(\aln\midd\bem)([\cbnm])_{[0]}\Bigr).
$$
$\phinm$\ will be the homomorphism mentioned in Definition-Theorem 0.1.
(The enumerative invariant for $g\geq 1$\ can be defined similarly
except that we know no interesting examples where
$\text{Mor}(\sig_g,B)$\ has the expected dimension, besides $g=1$\ and
$X=\bold P^n$.)

In the remainder of this section, we shall use the Poincare duality to
express the values of (1.7) as intersections of cycles in $\cbnm$.
Let $\al_i,\be_j\in A\sta X$\ be the cohomology classes as before.
Since $\phinm$\ is a homomorphism, it suffices to examine the value of
$\phinm$\ on those classes whose Poincare duals can be represented by
subvarieties of $X$.
Let $\Yn\sub X$\ and $\Zm\sub X$\ be subvarieties whose
Poincare duals are $\aln$\ and $\bem$\ respectively. (By this we mean
$\al_i$\ is the Poincare dual of the cycle $[Y_i]$\ and
likewise for $\be_j$.)
For $\zn\in\sig$, we consider the subscheme
$$Y=
\Bigl(\prod\indexn\{z_i\}\times Y_i\Bigr)
\times\Bigl(\prod\indexm Z_j\Bigr)\sub\sixnm
\tag 1.8
$$
and the subscheme
$$\Int\{\zdot,\Ydot,\Zdot\}=
\bl P^{n|m}\br\umo\bl Y\br
\sub\cbnm.
\tag1.9
$$
A moment of thought tells us that
closed points of this set are tuples
$$(C,x_1\cdotss x_{n+m}),\ \text{where}\ x_1\cdotss x_{n+m}\in C,
$$
such that $\ps(x_i)=z_i$\ and $\px(x_i)\in Y_i$\ for $1\leq i\leq n$\
and $\px(x_{n+j})\in Z_j$\ for $1\leq j\leq m$.

\pro{Definition 1.1}
We say the scheme $W_B$\ is good if $W_B$\ is a purely $\rho(B)$-dimensional
scheme and we say the collection $\{\zdot,\Ydot,\Zdot\}$\ is generic
(with respect to $W_B$) if in addition to $W_B$\ being good
there is a dense open subset
$$U\sub\Int\{\zdot,\Ydot,\Zdot\}
$$
such that
\roster
\item
$U$\ is contained in $\pi\umo(W_B\topo)$,
where $\pi\mh\cbnm\to W_B$\ is the projection;
\item
$U$\ has pure dimension $\rho(B)+m-(\sum \codim Y_i+\sum \codim Z_j)$\ and
\item
$U$\ is Cohen-Macaulay (smooth when $\char K=0$).
\endroster
\endpro

When $W_B$\ is good, $\{\zdot,\Ydot,\Zdot\}$\ is generic and
$$\dim \cbnm=\sum \codim Y_i+\sum \codim Z_j+n,
$$
then the following theorem states that
$\phinm(\aln\midd\bem)$\ is
the length (the number of points when $\char K=0$) of the finite scheme
$$\Int\{\zdot,\Ydot,\Zdot\}.
$$
(The length of a finite scheme is the length of its coordinate
ring as $K$-module.)

\pro{Theorem 1.2}
Let $\aln$\ and $\bem\in A\sta X$\ be cohomology classes whose
Poincare duals are represented by subvarieties $\Yn$\ and
$\Zm\sub X$\ respectively. Let $B\in\aotilalg$. Assume
$W_B$\ is good, the collection $\{\zdot,\Ydot,\Zdot\}$\
is generic and
$$\sum_{i=1}^n\codim Y_i+\sum_{j=1}^m\codim Z_j=\rho(B)+m,
$$
then
$$\phinm\left(\aln\mid\bem\right)=
\length\left(\Int\{\zdot,\Ydot,\Zdot\}\right).
$$
When $\char K=0$, then it is $\#\bl\Lam\br$,
where $\Lam$\ is as in Definition-Theorem 0.1.
\endpro

\proof
The first part follows directly from the formal properties of $A\sta X$\ and
the
Poincare duality [Fu,\S17.4]: Let $Y$\ be the subvariety in (1.8).
Since the Poincare dual $[Y]\dual$\ is
$$\Bl\prod\indexn e_0\time \al_i\Br\times\Bl\prod\indexm\be_j\Br
\in A\sta\Bl\sixnm\Br,
$$
$\vpsi_B(\aln\midd\bem)=\bl\Pnm\br\sta\bl[Y]\dual\br$.
By assumption, $\bl\Pnm\br\umo(Y)$\ is contained in
$\pi\umo(W_B\topo)$\ and has codimension equal to the
codimension of $Y$. Hence by [Fu,\S17.4 and \S8.3],
$$\bl\Pnm\br\sta\bigl([Y]\dual\bigr)\bigl[\cbnm\bigr]=
\bigl[\bl\Pnm\br\umo(Y)\bigr]\in A_0\cbnm
$$
and its degree is the length of $\bl\Pnm\br\umo(Y)$, because it is
Cohen-Macaulay. This proves the first part of theorem. Next we assume
$\bl\Pnm\br\umo(Y)$\ is smooth. Then its length is just the number
of points in it.
However, elements $(C,x_1\cdotss x_{n+m})$\ in $C_B^{n+m}$\
is in $\bl\Pnm\br\umo(Y)$\ if and only if $p_{\sig}(x_i)=z_i$\ and
$p_X(x_i)\in Y_i$\ for $i\leq n$\ and $p_X(x_{n+j})\in Z_j$\ for $j\leq m$.
Thus $\Int\{\zdot,\Ydot,\Zdot\}$\ is isomorphic to $\Lam$\ as
sets. Therefore,
the length of $\Int\{\zdot,\Ydot,\Zdot\}$\ is equal to $\#\bl\Lam\br$,
since it is smooth. This completes the proof of the theorem.
\qed

\pro{Corollary 1.3}
Let the notation be as in theorem 1.2 except that $K$\ is
the field of complex numbers. Then
$$\vphi_B\bl\aln\mid\bem\br=
\vphi_{[B]}^{\text{Sym}}\bl\mu(\al_1)\cdotss\mu(\al_n)\mid
\mu(\be_1)\cdotss\mu(\be_m)\br,
$$
where $\mu\mh A^pX\to H_{2p}(X)$\ is the map defined in
Theorem 0.3 and $\vphi_{[B]}^{\text{Sym}}$\ is the mixed
quantum invariants defined in [RT].
\endpro

\proof
Let $(X^{\text{top}},J_t)$, $t\in[0,1]$,
be a general (continuous) family of almost
complex structures satisfying $(X^{\text{top}},J_0)=X$.
We introduce
$$\Lam_t=\left\{(f,w_1\cdotss w_m)\Biggm\vert
\matrix f\mh\sig\to X^{\text{top}}\ \text{is $J_t$-holomorphic};\,
w_1\cdotss w_m\in\sig;\\
f(z_i)\in Y_i,\  1\leq i\leq n;\
f(w_j)\in Z_j,\  1\leq j\leq m
\endmatrix
\right\}.
$$
Since $\{\zdot,\Ydot,\Zdot\}$\ is generic, $\Lam_0=\Lam$\ is smooth.
It is proved in [RT] that
$\cup_{t\in[0,1]}\Lam_t$\ is a corbordism between $\Lam_0$\ and
$\Lam_1$. Therefore,
$$\vphi_{[B]}^{\text{Sym}}\bl\mu(\al_1)\cdotss\mu(\al_n)\mid
\mu(\be_1)\cdotss\mu(\be_m)\br=\deg \bl\Lam_1\br
$$
is exactly $\#\bl\Lam_0\br$, which is $\vphi_B\bl\aln\mid\bem\br$\
by Theorem 1.2.
This completes the proof of the corollary and Theorem 0.3.
\qed

\head 2. Rational curves in homogeneous manifolds
\endhead

In this section, we will collect some facts about rational
curves in homogeneous varieties that are essential to the
proof of the composition law.
Some of these results are standard and
can be found in literatures. For the convenience of the
readers, we will state them in accordance with our application and
will give reference or proof when necessary.

We continue to use the notion developed in section 1.
Namely, $X=G/P$\ is a smooth $d$-dimensional homogeneous variety,
$\sig$\ is the smooth rational curve, $B\in\aotilalg$\ and $\rho(B)
=c_1(X)\cdot B+d$.
As before, we let $H_B$\ be the Hilbert scheme of
curves $C\sub\sig\time X$\ satisfying (1.1) and let $\wbcirc\sub W_B
\sub H_B$\ be subsets defined after (1.2).
We still denote by $\ps$\ and $\px$\ the first and second
projection of $\sig\time X$\ and for
curve $C\sub\sig\time X$\ we denote by $\varphis$\ and
$\varphix$\ the induced morphism from $C$\ to $\sig$\ and $X$\
respectively.

We begin with a quick review of smoothness criterion of Hilbert
schemes.

\pro{Lemma 2.1}
Suppose $G$\ is a linear algebraic group over $K$\ and
$X=G/P$\ with $P\sub G$\ closed.
Then $\wbcirc$\ is a smooth quasi-projective
scheme of pure dimension $\rho(B)$.
\endpro

\proof
This is an easy consequence of deformation theory. Let $w\in\wbcirc$\ be
any point corresponding to the subscheme $C\sub \sig\time X$. Since $C$\
is irreducible, $C$\ is a smooth rational curve. Then according to [Gr, p21],
the the tangent space of $\wbcirc$\ at $w$\ is
$$\Hom\bl I_C/I_C^2,\OO_C\br,
\tag 2.1
$$
where $I_C$\ is the ideal sheaf of $C\sub \sig\time X$. Further $\wbcirc$\
is smooth at $w$\ if the obstruction
$$\Ext^1\bl I_C/I_C^2,\OO_C\br=0.
\tag 2.2
$$
Since $\varphi_{\sig}\mh C\to \sig$\ is an isomorphism,
$I_C\cong\varphi_X\sta\Omega_X$. On the other
hand, $\Omega_X\dual$\ is generated by global
sections since $G$\ is a linear algebraic group and $P\sub G$\ is closed
[Sp, 5.2.3]. Therefore, (2.2) must be trivial. Finally,
using Riemann-Roch theorem we calculate the dimension of (2.1)
to be $\rho(B)$, because of (2.2). This proves the lemma.
\qed

We now use generic position result to show that
for homogeneous variety the homomorphism $\vphi_B$\
can be defined as in Definition-Theorem 0.1.

\pro{Proposition 2.2}
Let $X$\ and $G$\ be as in Lemma 2.1 and let $B\in\aotilalg$. Assume
$\Yn$\ and $\Zm\sub X$\ are subvarieties satisfying
$$\sum_{i=1}^n\codim Y_i+\sum_{j=1}^m\codim Y_j=\rho(B)+m,
$$
then for general
$$(\zn,\sindh,\epm)\in\sig^n\time G^n\time G^m,
$$
the tuple $\{\zdot,\Ydot^{\si_{\cdot}},\Zdot^{\ep_{\cdot}}\}$\
is in generic position with respect to $W_B$\ as defined
in Definition-Theorem 0.1.
\endpro

The tool we need in studying translations of subvarieties
is the following general position result [Kl].

\pro{Lemma 2.3}
Suppose a connected algebraic group $H$\ acts transitively on a variety
$W$\ over $K$. Let $f\mh Y\to W$, $g\mh
Z\to W$\ be morphisms of varieties $Y$, $Z$\ to $W$. For each point
$\si$\ in $H$, let $Y^{\si}$\ denote $Y$\ with the morphism $\si\circ
f$\ from $Y$\ to $W$. Then
\roster
\item
for general $\si\in G\ucirc$,
$Y^{\si}\times_W Z$\ is either empty or of pure dimension
$$\dim(Y)+\dim(Z)-\dim(W),
$$
\item
if $Y$\ and $Z$\ are non-singular, and $\char(K)=0$\ (resp.
$\char(K)>0$), then for general $\si\in G\ucirc$, $Y^{\si}\times_W Z$\ is
non-singular (resp. Cohen-Macaulay).
\endroster
\endpro

Note that (1) of the lemma applies to any scheme $Y$\ and $Z$\
if we replace {\sl of pure dimension} by {\sl of dimension at most}.
Here for any scheme $Y$\ by dimension of $Y$\ we mean the maximum of the
dimensions of its irreducible components.

\noindent
{\sl Proof of Proposition 2.2}.
Consider the projection
$$P^{n|m}:C_B^{n+m}\lra \sixnm
$$
defined in (1.3) and the subvariety
$$\widetilde Y=\Bl\prod\indexn\{z_i\}\time \yisi\Br\time
\Bl\prod\indexm Z_j^{\ep_j}\Br\sub\sixnm.
$$
Since $\wbcirc\sub W_B$\ is dense, both have dimensions
$\rho(B)$, by lemma 2.1. Hence the dimension of $C_B^{n+m}$\ is the
same as the codimension of $\widetilde Y$\ in $\sixnm$,
by assumption on dimensions of $Y_i$'s and
$Z_j$'s. However, $\sixnm$\ is acted on transitively by
$\bl PGL(2)\time G\br^n\time G^m$.
Therefore, for general
$$(\zn,\sindh,\epm)\in \sig^n\time G^n\time G^m,
$$
the scheme
$$\Int\{\zdot,\Ydot^{\si_{\cdot}},\Zdot^{\ep_{\cdot}}\}=
\bl P^{n|m}\br\umo\bl\widetilde Y\br
$$
is discrete, is contained in $\pi\umo(\wbcirc)$\ and is
Cohen-Macaulay (smooth when $\char K=0$), by
repeatedly applying Lemma 2.3. This shows that the tuple
$\{\zdot,\Ydot^{\si_{\cdot}},\Zdot^{\ep_{\cdot}}\}$\ is in generic
position with respect to $W_B$.
\qed

\pro{Corollary 2.4}
Let $X$\ be as in Lemma 2.1. Then for any $B\in\aotilalg$\ the
enumerative invariant
$$\vphi_B(\cdot\midd\cdot):\prod\indexn A\sta X\times\prod\indexm
A\sta X\lra \ZZ
$$
can be defined as in Definition-Theorem 0.1
\endpro

Another smoothness criterion we need concerns the relative
Hilbert scheme $\wb$\
mentioned in the introduction. The exact formulation of this scheme
will be given in the next section. Here, we will prove a lemma that
will be useful in studying its smoothness.
Let $0\in V$\ be a smooth curve with uniformizing parameter $t$\
and let $Z$\ be a blowing-up of $\sig\time V$\
along a point in $\sig\time\{0\}$. $Z$\ is a flat family of curves over $V$.
We denote its fiber over $0\in V$\ by $Z_0$.

\pro{Lemma 2.5}
Let $X=G/P$\ be as in Lemma 2.1.
Assume $C\sub Z\time X$\ is a curve contained in $Z_0\time X$\
that is isomorphic to $Z_0$\ via projection $p_Z\mh Z\time X\to Z$,
\roster
\item
then $\Ext^1_C\bl I_C/I^2_C,\OO_C\br=0$, where
$I_C$\ the ideal sheaf of $C\sub Z\time X$;
\item
For $T_2\sub V$, where $T_2=\spec K[t]/(t^2)$,
there is a flat deformation $C_2$\ of $C$\ contained in $Z\time X$\
over $T_2$\ that makes the following diagram commutative and the
left square a Cartesian product:
$$\spreadmatrixlines{5pt}
\matrix C & \sub & C_2 & \sub & Z\time X\\
\downarrow && \downarrow && \downarrow \\
0 & \in & T_2 & \sub & V.\\
\endmatrix
$$
\endroster
\endpro

\proof
Since $\varphi_{Z_0}\mh C\to Z_0$\ induced by $p_Z$\ is an isomorphism,
$C$\ is the graph of $p_X\circ\varphi_{Z_0}\umo\mh Z_0\to X$. Hence
if we let ${}'\!I_C$\ be the ideal sheaf of $C\sub Z_0\time X$, then
${}'\!I_C$\ is isomorphic to $\varphi_X\sta\Omega_X\dual$\ and
$I_C/I_C^2$\ belongs to the exact sequence
$$0\lra \OO_C\lra I_C/I_C^2\lra {}'\!I_C/'\!I_C^2\lra 0.
$$
Hence (1) follows because $\Ext^1_C('\!I_C/'\!I_C^2,\OO_C)$\ and
$\Ext^1_C(\OO_C,\OO_C)$\ are trivial.

Next we prove (2). Because $C$\ is a graph of $Z_0\to X$, $C$\ is a
local complete intersection. Hence locally $C\sub Z\time X$\ can
be extended to curves in $Z\time X$\ flat over $T_2$.
By deformation theory (see for example [Gr, p21] or [Ko, \S 1.2]),
the obstruction to the existence of such $C_2$\ lies in
$\Ext^1_C('\!I_C/'\!I_C^2,\OO_C)$, which is trivial by the previous
argument. This prove the second part of the lemma.
\qed

We now state and prove our main technical lemma
concerning the set of curves in $Z_0\time X$\
that intersect with a set of prescribed subvarieties.
Here $Z_0$\ as before is a nodal curve that is the union of two
rational curves $\sigo$\ and $\sigt$\ intersecting along $x_0\in\sigo$\
and $y_0\in\sigt$. We fix subvarieties $\Yn$\ and $\Zm\sub X$.
For integer $k$\ and $l\geq 1$, we define $\r_k^l$\ to be the set of all
reduced,
connected curves $C\sub Z_0\time X$ such that in addition to $C$\ being
unions of exactly $l$\ rational curves,
$$\deg\bl C,p_{Z_0}\sta\tilde\OO_{\sigo}(1)\br=
\deg\bl C,p_{Z_0}\sta\tilde\OO_{\sigt}(1)\br=1,\
\deg\bl C,\px\sta(-K_X)\br=k.
\tag 2.3
$$
Here $\tilde\OO_{\sigo}(1)$\ is the invertible sheaf on $Z_0$\
having degree 1 along $\sigo$\ and 0 along $\sigt$, likewise for
$\tilde\OO_{\sigt}(1)$. Further, for fixed partition $n=n_1+n_2$, $n_1,n_2\geq
0$,
and $x=(x_1\cdotss x_{n_1})\in\sigo^{n_1}$,
$y=(y_1\cdotss y_{n_2})\in\sigt^{n_2}$, $\si=(\sindh)\in G^n$\ and
$\ep=(\ep_1\cdotss \ep_m)\in G^m$, we define
$$\r_k^l(\si,\ep)=\left\{C\in\r_k^l\Bigm\vert
\matrix
C\cap\{x_i\}\time Y_i^{\si_i}\ne\emptyset\ \forall 1\leq i\leq n_1;\
C\cap\{y_j\}\time Y_{n_1+j}^{\si_{n_1+j}}\ne\emptyset\\
\forall 1\leq j\leq n_2;\
C\cap (Z_0\time Z_h^{\ep_h})\ne\emptyset\ \forall 1\leq h\leq m.
\endmatrix
\right\}.
$$

\pro{Proposition 2.6}
With the notation as before. Then for general $(x,y,\si,\ep)\in
\sigo^{n_1}\time \sigt^{n_2}\time G^n\time G^m$,
the set $\r_k^l(\si,\eps)$\ satisfies
$$\dim \r_k^l(\si,\eps)\leq k+d+m-(l-2)-\sum_{i=1}^n\codim_X Y_i-
\sum_{j=1}^m\codim_XZ_j.
$$
\endpro

\pro{Corollary 2.7}
Let $X=G/P$\ and $\Yn,\Zm\sub X$\ be as before such that
$$\sum_{i=1}^n\codim Y_i+\sum_{j=1}^m\codim Z_j=k+d+m,
$$
where $k\geq 1$.
Then for any partition $n=n_1+n_2$, there is a dense open subset
$$U\sub \sigo^{n_1}\time \sigt^{n_2}\time G^n\time G^m
$$
such that for any $(x,y,\si,\ep)\in U$\
\roster
\item
$\r_{k_0}^l(\si,\ep)=\emptyset$\ for all $k_0<k$\ and $l\geq 1$;
\item
$\r_k^l(\si,\ep)=\emptyset$\ for all $l\geq 3$;
\item
$\r_k^2(\si,\ep)$\ is a finite set and further for any element
$C\in R^2_k(\si,\ep)$, $C\cap\bl Z_0\time Z_j^{\ep_j}\br$\ is
contained in the smooth locus of $C$\ and
$$C\cap\bl Z_0\time (Z_{j_1}^{\ep_{j_1}}\cap Z_{j_2}^{\ep_{j_2}})\br=
\emptyset
$$
for all $j_1\ne j_2$.
\endroster
\endpro

Before we prove this proposition, we need to study curves in $\sig\time X$.
We first consider the case where
$W$\ is a variety acted on transitively by a linear algebraic
group $H$\ and $S_1\cdotss S_h$\ are subvarieties of $W$.
Let $\s^l_k$\ be the set of degree $k$\ (calculated using $c_1(W)$)
reduced connected curves $C\sub X$\ that are unions of exactly $l$\
rational curves and let
$$\s_k^l\bl S_i\mid i\in\Lam\br=
\bigl\{C\in \s^l_k\mid C\cap S_i\ne\emptyset\ \text{for}\ i\in\Lam\bigr\},
\tag 2.4
$$
where $\Lam$\ is the set $\{1\cdotss h\}$.

\pro{Lemma 2.8}
Let $W$\ and $S_i$\ be as before. Then for general
$\ups=(\upsilon_1\cdotss\upsilon_h)\in H^h$,
$$\dim \s_k^l\bl \siui\mid i\in\Lam\br\leq
k+\dim W-(l+2)-\sum_{i\in\Lam}(\codim S_i-1).
$$
\endpro

\proof
We will prove the lemma by induction on $l$. We begin with
$\s_k^1(S_i\mid i\in\Lam)$. Let
$$Z=\{(x,C)\mid C\in \s^1_k,\ x\in C^h\}
$$
be endowed with the obvious reduced scheme structure and let
$f\mh Z\to W^h$\ be the morphism that sends $(x,C)$\ to $x\in W^h$\ via
$C^h\sub W^h$. Also, we let $Y=\prod_{i\in\Lam} S_i$\ and let
$g\mh Y\to W^h$\ the obvious inclusion. Note that $H^h$\ acts
transitively on $W^h$.
For any $\ups=(\ups_1\cdotss \ups_h)\in H^h$, any elements
$$\bl(x,C),y\br\in Z\times_{W^h}Y^{\ups}
$$
satisfies $x_i=\ups_i(y_i)$\ for all $i\in\Lam$, by the
definition of Cartesian product. Hence $C\cap
\ups_i(S_i)\ne\emptyset$\ and then $C\in\s_k^1(S_i^{\ups_i}\mid i\in\Lam)$.
Now let $P\mh Z\times_{W^h}Y^{\ups}\to \s^1_k$\
be the map sending $\bl(x,C),y\br$\ to $C$. It is easy to see that
$$\s^1_k\bigl(S_i^{\ups_i}\mid i\in\Lam\bigr)=P\bl Z\times_{W^h}Y^{\ups}\br.
$$
Finally, after applying Lemma 2.3 to morphism $f$\ and $g$, we
conclude that for general $\ups\in H^h$,
the set $Z\time_{W^h} Y^{\ups}$\ and then
$\s_k^1(S_i^{\ups}\mid i\in\Lam)$\ have dimensions at most
$$\dim Z+\sum_{i\in\Lam}\dim S_i -h\dim W =k+\dim W-3-\sum_{i\in\Lam}
(\codim S_i-1).
$$
Here we have used the fact that $\dim \s_k^1\leq k+\dim W-3$\ that follows
from Lemma 2.1. This proves the lemma for $l=1$.

Next, we prove the lemma for $l>1$. Let $1\leq k_1<k$\ be any
integer and let $\Lam_1\sub\Lam$\ be any subset.
For $\ups=(\ups_1\cdotss\ups_h)\in H^h$, we consider the subset
$$\tkls\sub\s_{k_1}^{l-1}(\siui\mid i\in\Lam_1)\times
\s_{k_2}^{1}(\siui\mid i\in\Lam_2)
$$
that consists of
pairs of curves $(C_1,C_2)$\ such that $C_1\cap C_2\ne\emptyset$\
and let
$$\tkls^0\sub\tkls
$$
be the subset of those $(C_1,C_2)$\ such that $C_1$\ and $C_2$\
share no common components.
Here $k_2=k-k_1$, $\Lam_2=\Lam-\Lam_1$\ and $\s_{k_1}^{l-1}(\siui
\mid i\in\Lam_1)$\ and $\s_{k_2}^{1}(\siui\mid i\in\Lam_2)$\ are
defined as in (2.4). Let $\f$\ be the map
$$\f:\tkls^0\lra \s_k^l\bl\siui\mid i\in\Lam\br
$$
sending $(C_1,C_2)\in\tkls^0$\ to $C_1\cup C_2$. Clearly
$$\opcup_{1\leq k_1<k,\ \Lam_1\sub\Lam} \f\bl\tkls^0 \br
=\s_k^l\bl\siui\mid i\in\Lam\br,
$$
since every element $C\in\s_k^l(S_i^{\ups_i}\mid i\in\Lam)$\ splits into
$C_1\cup C_2$\ with $C_1$\ irreducible and $C_2$\ connected.
Hence to prove the lemma it suffices to show that for any choice of
$k_1$\ and $\Lam_1$,
$$\dim\tkls\leq k+\dim W-(l+2)-\sum_{i\in\Lam}(\codim Y_i-1)
\tag 2.5
$$
holds for general $\ups\in H^h$.

Now we prove (2.5). Without loss of generality, we can assume
$\Lam_1=\{1\cdotss r\}$\ and then $\Lam_2=\{r+1\cdotss h\}$.
By induction hypothesis, we can find a non-empty open subset $U_1\sub H^r$\
such that for $\upsilon\pri=(\upsilon_1\cdotss\upsilon_r)\in U_1$,
$$\dim \s_{k_1}^{l-1}\!\bl\siui\mid i\leq r\br\leq
k_1+\dim W-\bl (l-1)+2\br-\sum_{i=1}^r(\codim S_i-1).
\tag 2.6
$$
For similar reason, we can find a non-empty open subset $U_2\sub H^{n-r}$\
such that for any $\upsilon\doublepri=(\ups_{r+1}\cdotss\ups_n)
\in U_2$,
$$\dim \s_{k_2}^1\!\bl\siui\mid i\geq r+1\br\leq
k_2+\dim W-\bl 1+2\br-\sum_{i=r+1}^h(\codim S_i-1).
\tag 2.7
$$
As before, we consider the set
$$\z_1=\bigl\{(\ups\pri,z_1,C_1)\mid \ups\pri\in U_1,\
C_1\in\s_{k_1}^{l-1}\!\bl\siui\mid i\leq r\br\
\text{and}\ z_1\in C_1\bigr\}
$$
and $f_1\mh \z_1\to W$\ that is the map sending $(\ups\pri,z_1,C_1)$\
to $z_1$. Similarly, we let $\z_2$\ be the set
$$\z_2=\bigl\{(\ups\doublepri,z_2,C_2)\mid
\ups\doublepri\in U_2,\
C_2\in\s_{k_2}^{1}\bl\siui\mid i\geq r+1\br\
\text{and}\ z_2\in C_2\bigr\}
$$
and let $f_2\mh \z_2\to X$\ be the map sending $(\ups\doublepri,z_2,C_2)$\
to $z_2$. Since $\z_1$\ and $\z_2$\
are finite unions of varieties (both endowed with reduced
scheme structures) and $f_1$\ and $f_2$\ are morphisms,
we can apply Lemma 2.3 to conclude that there is a non-empty open subset
$H^{\circ}\sub H$\ such that for $\tau\in H^{\circ}$,
$$\dim\z_1^{\tau}\time_W\z_2\leq
\dim\z_1+\dim\z_2-\dim W.
$$
Combined with (2.6) and (2.7), we obtain
$$\dim \z_1^{\tau}\time_W\z_2\leq
k+\dim W+h\dim H-(l+2)-\sum_{i\in\Lam}(\codim S_i-1).
\tag 2.8
$$
Finally, let
$$\h: \z_1^{\tau}\time_W\z_2
\lra \s_{k_1}^{l-1}\bl \siui\mid i\leq r\br\times
\s_{k_2}^{1}\bl \siui\mid i\geq r+1\br\times H^h
$$
be the map sending
$$\bigl\{(\ups\pri,z_1,C_1),(\ups\doublepri,z_2,C_2)\bigr\}
\in \z_1^{\tau}\times_W\z_2
$$
to $\bl \tau(C_1),C_2,(\ups\pri,\ups\doublepri)\br$. Since
$\tau(z_1)=z_2$, $\tau(C_1)\cap C_2\ne\emptyset$\ and vice versa
if $C_1\cap C_2\ne\emptyset$, then there are $z_1\in \tau\umo(C_1)$\
and $z_2\in C_2$\ such that $\tau(z_1)=z_2$. Hence the image of
$\h$\ is exactly
$$\opcup_{\ups\in U_1^{\tau}\time U_2}
\Gamma\bl k_1, \Lam_1, \eps\br\time\{\ups\},
$$
where $U_1^{\tau}$\ is the translation of $U_1$\ under $\tau$\
acting diagonally on $H^r$. Therefore
$$\dim \opcup_{\ups\in U_1^{\tau}\time U_2}
\Gamma\bl k_1, \Lam_1, \ups\br\time\{\ups\}\leq
k+\dim W+h\dim H-(l+2)-\sum_{i\in\Lam}(\codim S_i-1),
$$
by (2.8). However, since $\cup_{\tau\in H^{\circ}}U_1^{\tau}\time U_2$\
is dense in $H^h$,
there is an open subset $U\sub H^h$\ such that for any $\ups\in U$,
$$\dim \Gamma\bl k_1, \Lam_1, \ups\br \leq
k+\dim W-(l+2)-\sum_{i\in\Lam}(\codim S_i-1)
$$
as desired. This completes the proof of the lemma.
\qed

In applying Lemma 2.8 to prove Proposition 2.6, we will consider the
case $W=\Sigma\times X$\ and the set $\t^l_k$\ that consists of curves
$C$\ that are reduced, connected and are unions of exactly $l$\ rational
curves such that
$$\deg\bl\ps\sta\OO_{\Sigma}(1), C\br=1\quad \text{and}\quad
\deg\bl\px\sta(-K_X), C\br=k.
$$

\noindent
{\sl Proof of Proposition 2.6}.
We first introduce more sets. Let $0\leq k_1\leq k$, $l_1\geq 1$\ be
integers and $\Lam_1\sub\Lam$\ be any subset, where $\Lam=\{1,\cdotss
m\}$\ this time. For $\si\in G^n$\ and $\eps\in G^m$, we define
$\t_{k_1}^{l_1}(\si,\eps,\Lam_1)\sub\t_{k_1}^{l_1}$\ to be the subset
of curves $C\in \t_{k_1}^{l_1}$\ such that
$$C\cap(\{x_i\}\times\yisi)\ne\emptyset\  \forall\, 1\leq i\leq n_1\
\text{and}\  C\cap(\sigo\times\zjej)\ne\emptyset\  \forall\, j\in\Lam_1.
\tag 2.9
$$
We let $k_2=k-k_1$\ and $\Lam_2=\Lam-\Lam_1$. For any $l_2\geq 1$,
we define $\t_{k_2}^{l_2}(\si,\eps,\Lam_2)\sub\t_{k_2}^{l_2}(\si,\ep,\Lam_2)$\
to be the subset of curves $C\in\t_{k_2}^{l_2}$\ such that $C$\
satisfy (2.9) with $\{x_i\}\time\yisi$\ (resp. $n_1$, resp. $\Lam_1$)
replaced by $\{y_i\}\time Y_{n_1+i}^{\si_{n_1+i}}$\
(resp. $n_2$; resp. $\Lam_2$).
Similar to the proof of Lemma 2.7, we consider
$$\Gamma(k_1,l_1,l_2,\si,\eps,\Lam_1)
\sub\t_{k_1}^{l_1}(\si,\eps,\Lam_1)\times\t_{k_2}^{l_2}(\si,\eps,\Lam_2)
$$
that is the set of pairs $(C_1,C_2)$\ such that $\tilde C_1\cap\tilde C_2\ne
\emptyset$, where $\tilde C_i\sub Z_0\time X$\ is the image of $C_i
\sub\sig\time X$\ under $\sig\time X=\sig_i\time X\sub Z_0\time X$.
(Here we fix an isomorphism $\sig=\sig_i$\ once and for all.) We let
$$\Gamma(k_1,l_1,l_2,\si,\eps,\Lam_1)^0\sub
\Gamma(k_1,l_1,l_2,\si,\eps,\Lam_1)
$$
be the subset consisting of pairs $(C_1,C_2)$\ such that
$\tilde C_1$\ and $\tilde C_2$\
have no common components. Because of (2.3),
every element $C\in\r_k^l(\si,\eps)$\ splits into two
connected curves $C_1\sub\sigo\time X$\ and $C_2\sub\sigt\time X$\
whose number of components add up to $l$.
Thus the map $\h$\ that sends $(C_1,C_2)$\ to $\tilde C_1\cup\tilde C_2$\
satisfies
$$\opcup_{\Sb l_1+l_2=l,\, 0\leq k_1\leq k\\ \Lam_1\sub\Lam\endSb}
\h\bl \Gamma(k_1,l_1,l_2,\si,\eps,\Lam_1)^0 \br
=\r_k^l(\si,\eps).
$$
Therefore, to prove the lemma it suffices to show that for general
$(x,y,\si,\ep)\in\sigo^{n_1}\time\sigt^{n_2}\time G^n\time G^m$,
$$\dim \Gamma(k_1,l_1,l_2,\si,\eps,\Lam_1)
\leq k+d+m-(l-2)-\sum_{i=1}^n\codim Y_i-\sum_{j=1}^m\codim Z_j.
\tag 2.10
$$

The proof of this inequality
is again based on Lemma 2.3 as we did in the proof of Lemma 2.8.
First, by Lemma 2.3, there is a non-empty open subset
$U\sub \sigo^{n_1}\time
\sigt^{n_2}\time G^n\time G^m$\ such that for any $(x,y,\si,\ep)\in U$,
$$\align
\dim\t_{k_1}^{l_1}(\si,\eps,\Lam_1)\leq&
k_1+d+1-l_1-\sum_{i=1}^{n_1}\codim Y_i
-\sum_{j\in\Lam_1}\bl\codim Z_j-1\br\\
\dim \t_{k_2}^{l_2}(\si,\ep,\Lam_2)\leq&
k_2+d+1-l_2-
\sum_{i=1}^{n_2}\codim Y_{n_1+i}
-\sum_{j\in\Lam_2}\bl\codim Z_j-1\br.
\endalign
$$
(Note that element in $\t^l_k$\ have degree $2+k$\ according to the
convention in Lemma 2.8.)
Now, let $\z_1$\ and $\z_2$\ be the sets defined as
$$\align
\z_1=\bigl\{(z_1,C_1)\mid &C_1\in
\t_{k_1}^{l_1}(\si,\eps,\Lam_1),\, z_1\in C_1\cap\{x_0\}\time X\bigr\},\\
\z_2=\bigl\{(z_2,C_2)\mid &C_2\in
\t_{k_2}^{l_2}(\si,\eps,\Lam_2),\, z_2\in C_2\cap\{y_0\}\time X\bigr\},\\
\endalign
$$
and $f_i$, $i=1,2$, be maps from $\z_i$\ to $X$\ sending $(z_i,C_i)$\
to $z_i$. $f_1$\ and $f_2$\ are morphisms after we endow $\z_1$\ and
$\z_2$\ with the reduced scheme structures. Note that the map
$$\z_1\time \z_2\to \Gamma(k_1,l_1,l_2,\si,\ep,\Lam_1)
$$
sending $\{(z_1,C_1),(z_2,C_2)\}$\ to $(C_1,C_2)$\
is surjective. Thus to establish (2.10) it suffices to show
$$\dim\z_1\times_X\z_2\leq k+d+m-(l-2)-\sum_{i=1}^{n}\codim Y_i
-\sum_{j=1}^m\codim Z_j.
\tag 2.11
$$
We will prove this by using Lemma 2.3.

Since $\sig$\ is homogeneous, possibly after shrinking $U$\ if
necessary, we can assume
$$\dim\z_1=\dim \t_{k_1}^{l_1}(\si,\eps,\Lam_1)\
\text{and}\ \dim\z_2=\dim \t_{k_2}^{l_2}(\si,\eps,\Lam_2).
$$
However, since $X$\ is homogeneous there is a non-empty
open set $U_0\sub U\times G$\
such that for any $(x,y,\si,\ep,\tau)\in U_0$,
$$\dim \z_1^{\tau}\times_X\z_2\leq \dim\z_1+\dim\z_2-\dim X,
$$
by Lemma 2.3. Put them together, we get
$$\dim\z_1^{\tau}\times _X\z_2\leq k+d+m-(l-2)-
\sum_{i=1}^{n}\codim Y_i -\sum_{j=1}^m\codim Z_j,
\tag 2.12
$$
since $l_1+l_2=l$. We now show (2.12) leads to (2.10). The argument
is a repetition of that of Lemma 2.4 by exploiting the homogeneity of
$X$. One can show that $\z_1^{\tau}\time_X\z_2$\ is canonically
isomorphic to $\z_1\pri\time_X\z_2$, where $\z_1\pri$\ is defined as
that of $\z_1$\ with
$$\si_i\ \text{replaced by}\ \tau\cdot\si_i\ \text{for}\ 1\leq
i\leq n_1\ \text{and}\ \ep_j\ \text{replaced by}\ \tau\cdot\ep_j\
\text{for}\ j\in\Lam_1.
\tag 2.13
$$
Since (2.12) hold for all $(x,y,\si,\ep,\tau)\in U_0$,
(2.11) will hold for $(x,y,\si\pri,\ep\pri)\in
\sigo^{n_1}\time \sigt^{n_2}\time G^n\time G^m$, where $\si\pri$\ and
$\ep\pri$\ are derived from $(x,y,\si,\ep,\tau)$\ by the rule (2.13).
However, all such $(x,y,\si\pri,\ep\pri)$\ still form a dense open
subset of $\sigo^{n_1}\time \sigt^{n_2}\time G^n\time G^m$.
Therefore, (2.11) holds for general $(x,y,\si\pri,\ep\pri)$\ and the
proposition follows.
\qed

\head 3. Composition law
\endhead

The composition law for $\phinm$\
runs as follows: Suppose the Poincare dual of the diagonal
$\Del \sub X\times X$\ admits a Kunneth decomposition
$$[\Del]\dual=\sum_{l=1}^p\zeta_l\times\tilde\zeta_l
\tag 3.1
$$
then for any partition $n=n_1+n_2$\ with $n_1,n_2\geq 2$,
$$\align
&\phinm(\aln\midd\bem)\\
=&\sum_{\Sb B=B_1+B_2,\,\tau\in S(m)\\
0\leq k\leq m,\,1\leq l\leq p\endSb}
    {1\over k!(m-k)!}\cdot
 \vphi_{B_1}\Bigl(\zeta_l,\al_1\cdotss \al_{n_1}\midd
      \be_{\tau(1)}\cdotss\be_{\tau(k)}\Bigr)\\
    &\qquad\qquad\qquad\times
\vphi_{B_2}\left(\tilde\zeta_l,\al_{n_1+1}\cdotss \al_n
\midd \be_{\tau(k+1)}\cdotss \be_{\tau(m)}\right).
\tag 3.2
\endalign
$$
The first mathematical proof of this formula for smooth Fano manifolds
is due to [RT], where they proved this relation (for all $g\geq 0$)\ by
using analytic method. The algebraic approach for rational surfaces and
Grassmannian appeared recently in [CM, Be].
In this section, we will prove the relation (3.2) (for $g=0$) for all
homogeneous varieties over $K$\ by degeneration methods
in algebraic geometry. (Note that by work of [FMS$^2$], the Kunneth
decomposition always exists for the varieties considered in Theorem 0.2.)

We keep the notation developed in the previous sections. Namely,
$X=G/P$\ is a dimension $d$\ smooth homogeneous variety over $K$\ for a
linear algebraic group $G$\ and a closed subgroup $P\sub G$. Let
$\Yn$\ and $\Zm$\ be subvarieties of $X$\
and $\aln$\ and $\bem\in A\sta X$\
be the Poincare duals of the cycles represented by these
subvarieties. We denote $e_0\in A^1\sig$\
the poincare dual to a point in $\sig$. We agree that by a
subvariety $Y\sub X$\ in general position we mean $Y$\ is a
general member among all translations $\{Y^{\si}\midd\si\in G\}$.
Let $B\in\aotilalg$.

We first introduce the degeneration of the Hilbert scheme $W_B$.
Let $V$\ be a rational curve as a parameter
space and let $0\in V$\ be fixed.
We form a family of curves $Z$\ over $V$\ by blowing up a
point in $\stimev$\ lying over $0\in V$.
Let $\pi\mh Z\to V$\ be the projection and let $Z_v=\pi\umo(v)$\
be the fiber over $v\in V$. $Z_v$\ is isomorphic to $\sig$\
for $v\ne0$\ and $Z_0$\
is a union of two rational curves, which we denote by $\sigo$\ and
$\sigt$. Let $x_0\in\sigo$\ (resp. $y_0\in\sigt$) be the point
corresponding to the singular point in $Z_0$.
We also fix a projection $h_1\mh Z\to \sig$\ that is the result of
blowing down $Z$\ along $\sigt$\ and then project to $\sig$. We fix
an $h_2\mh Z\to\sig$\ similarly by blowing down $\sigo$\ first.
We now define a relative functor
$$\Cal F:\bigl\{\text{category of schemes}/V\bigr\}\lra
\bigr\{\text{Sets}\bigr\}/\sim
$$
that sends any scheme $T$\ over $V$\ to the set of all subschemes
$S\sub T\time_V Z\time X$\ flat over $T$\ of which the following
holds: For any closed $t\in T$\ over $v\in V$,
$$[\px(S_t)]\in B,\ \xx(\OO_{S_t})=1,\ \deg(h_1\sta\OO_{\sig}(1), S_t)=
\deg(h_2\sta\OO_{\sig}(1), S_t)=1.
\tag 3.3
$$
$\Cal F$\ is represented by the Hilbert scheme $\h_B$, projective
over $V$. Since $Z\to V$\ is a constant family over $V-0$,
for $v\ne 0$\ the fiber of $\h_B$\ over $v$\ is isomorphic
to $H_B$. We let $\wb\sub\h_B$\ be the closure of $W_B\times(V-0)
\sub\h_B$. $\wb$\ is projective and flat over $V$.

In the following, we will define the relative analogue of (1.5).
Let $\c_B\sub Z\time X\times\wb$\ be the (restriction of the) universal family
and let $\tcbnm$\ be the product of $n+m$\ copies of $\c_B$\ over $\wb$.
As before, we let $\piiz$\ and $\piix$\ be projections from $\tcbnm$\
to $Z$\ and $X$\ respectively that factor through the i-th copy of
$\tcbnm$. After choosing a partition $n=n_1+n_2$\ with $n_1,n_2\geq 2$,
we can define a morphism
$$P^{n_1:n_2\vert m}:\tcbnm\lra
\bl\sig\time X\br^{n_1}\time\bl\sig\time X\br^{n_2}\time X^m
$$
that is the product
$$\Bigl(\prod_{i=1}^{n_1}( h_1\circ\piiz) \time \piix \Bigr)\times
   \Bigl(\prod_{i=n_1+1}^{n}( h_2\circ\piiz) \time \piix \Bigr)\times
\Bigl(\prod_{j=1}^{m}\pi_{n+jX}\Bigr).
$$
We let
$$\prod\indexn A\sta X\times\prod\indexm A\sta X\lra
A\sta\Bl\sixnotm\Br
\tag 3.4
$$
be the inclusion sending $(\aln;\bem)$\ to
$$\Bl\prod_{i=1}^{n_1} e_0\time\al_i \Br\time
\Bl\prod_{i=n_1+1}^n e_0\time\al_{i} \Br\time
\Bl\prod\indexm\be_j\Br
\in A\sta\Bl\sixnotm\Br.
\tag 3.5
$$
($e_0$\ is the Poincare dual of a point.)
The composition of $\bl P^{n_1:n_2\vert m}\br\sta$\ with the
inclusion (3.4) defines a homomorphism
$$\tilde\vpsi_B^{n_1:n_2}:\prod\indexn A\sta X\times
\prod\indexm A\sta X  \lra A\sta\tcbnm
$$
that is the relative version of $\psinm$\ in (1.5). By definition of
cohomology theory, $\tilde\vpsi_B^{n_1:n_2}$\ paired with the
fundamental class $[\c^{n+m}_{B,v}]\in A\lsta\tcbnm$,
where $\c_{B,v}^{n+m}$\ is the fiber of $\tcbnm$\ over $v\in V$,
defines a homomorphism
$$\prod\indexn A\sta X\times \prod\indexm A\sta X
\lra A\lsta\tcbnm
$$
that composed with the degree homomorphism
$$A\lsta \tcbnm\mapright{\text{pr}} A_0\tcbnm\mapright{\deg} \ZZ
$$
gives rise to the  relative version of (1.7):
$$\tilde\vphi^{n_1:n_2}_{B,v}:
\prod\indexn A\sta X\times \prod\indexm A\sta X\lra \ZZ.
$$
In explicit form,
$$\tilde\vphi^{n_1:n_2}_{B,v}(\aln\midd\bem)=
\deg\Bl \tilde\vpsi_B^{n_1:n_2}(\aln\midd\bem)[\c_{B,v}^{n+m}]\Br_{[0]}.
$$
On the other hand, since $\c_{B}^{n+m}$\ is constant over $V-0$,
$$\align
\vphi_B(\aln\midd\bem)
=&\deg\Bl \vpsi_B(\aln\midd\bem)[C_{B}^{n+m}]\Br_{[0]}\\
=&\deg\Bl \tilde\vpsi_B^{n_1:n_2}(\aln\midd\bem)[\c_{B,v}^{n+m}]\Br_{[0]}
\endalign
$$
for $v\in V-0$. However,
$[\c_{B,v}^{n+m}]=[\c_{B,0}^{n+m}]$\ in $A\lsta \c_{B}^{n+m}$\
by our choice of $V$\ (it is rational) and the flatness of
$\c_{B}^{n+m}$\ over $V$. Hence
$$
\vphi_B(\aln\midd\bem)=\deg\Bl
\tilde\vpsi_B^{n_1:n_2}(\aln\midd\bem)[\c_{B,0}^{n+m}]\Br_{[0]}.
\tag 3.6
$$
Therefore in order to establish the
composition law, we only need to show that
the right hand side of (3.2) is identical to
$$\deg \Bl\tilde\vpsi_B^{n_1:n_2}(\aln\midd\bem)
[\c_{B,0}^{n+m}]\Br_{[0]}.
\tag 3.7
$$
We remark that we only need to consider the case when
$$\sum\indexn \codim Y_i+\sum\indexm \codim Z_j=\rho(B)+m,
\tag 3.8
$$
since otherwise both sides of (3.2) vanish automatically for dimension
reason. In the following, we will assume without further mentioning that
the identity (3.8) always hold.

We now construct explicitly the cycle in (3.7) using
Poincare duality, similar to that in the proof of
theorem 1.2.
We choose $\zn\in\sig$\ and consider the subvariety
$$\widetilde{Y}=\Bl\prod_{i=1}^{n_1} \{z_i\}\time Y_i\Br\times
\Bl\prod_{i=n_1+1}^{n}\{z_i\}\time Y_i\Br\times
\Bl\prod_{j=1}^m Z_j\Br \sub \sixnotm.
$$
Note that $\widetilde{Y}$\ represents the Poincare dual of (3.5).
By our assumption on (3.8), the codimension of
$\widetilde{Y}$\ is one less than the dimension of $\tcbnm$. In a moment,
we will show that when $z_i$'s, $Y_i$'s and $Z_j$'s are in
general positions, then the subscheme
$$\bl P^{n_1:n_2|m}\br\umo
\bl\widetilde{Y}\br
\sub\c_{B}^{n+m},
$$
is one-dimensional and is flat and Cohen-Macaulay over $V$\ at $0\in V$.
Hence (3.7) will be the same as the length of
$$\tIntnotm\bl\zdot,\Ydot,\Zdot\br=
\bl P^{n_1:n_2|m}\br\umo
\bl\widetilde{Y}\br
\bigcap \c_{B,0}^{n+m},
\tag 3.9
$$
since it is 0-dimensional and Cohen-Macaulay [Fu, Proposition 7.1
and Corollary 17.4].
A moment of thought tells us that points in (3.9) are curves
$S\sub Z_0\time X$\ whose two irreducible components as
pairs of curves lie in $C_{B_1}\udot\time C_{B_2}\udot$\
for some $B_1+B_2=B$\ (see Corollary 2.7).
It is the understanding of this
correspondence that will lead to the proof of the composition
law.

We first relate the summands in the right hand side of (3.2)
to the lengths of subschemes in $C_{B_1}\udot\time C_{B_2}\udot$. To
this end, we introduce some convention that will be used
in the remainder of this section. Let
$B=B_1+B_2$, $n=n_1+n_2$\ ($n_1,n_2\geq 2$\
as always) and $m=m_1+m_2$\ ($m_1,m_2\geq 0$) be partitions
that will be
fixed momentarily. For convenience, we will think of $C_{B_1}$\
as a universal family in $\sigo\time X\time W_{B_1}$\ (namely points in
$C_{B_1}$\ are pointed curves in $\sigox$) and
$$C_{B_1}^{1+n_1+m_1}=C_{B_1}\time_{W_{B_1}}\cbnmo
$$
a product of $(1+n_1+m_1)$\ copies of $C_{B_1}$\ indexed by $i=0$\
(corresponding to the first copy of $\cbo$) and $i=1\cdotss n_1+m_1$\
(corresponding to $n_1+m_1$\ copies in $\cbnmo$). We adopt the
same convention to $\cbnmto$.
Under this convention, we have the obvious projections
$$\piso:\cbnmoo\lra\sigo\ \text{and}\ \pix:\cbnmoo\lra X
$$
indexed from $0$\ to $n_1+m_1$\ (that factor through the $i$-th copy of
$\cbnmoo$) and
$$\qist:\cbnmto\lra\sigt\ \text{and}\ \qix:
\cbnmto\lra X
$$
alike indexed from $0$\ to $n_2+m_2$. Let $\zn\in\sig$\ be
points away from $h_1(\sigt)$\ and $h_2(\sigo)$, let $x_i=h_1\umo(z_i)
\in\sigo$\ for $i\leq n_1$\ and let $y_i=h_2\umo(z_{n_1+i})\in\sigt$\
for $i\leq n_2$.
Similar to (1.3), we introduce projection
$$P^{n_i|m_i}:C_{B_i}^{1+n_i+m_i} \lra
\bl\sig_i\time X\br^{n_i}\time X^{m_i},\quad i=1,2
\tag 3.10
$$
that factor through $C_{B_i}\time_{W_{B_i}}C_{B_i}^{n_i+m_i}
\to C_{B_i}^{n_i+m_i}$\
and subvarieties
$$\align
\Theta_1(\tau)&=\Bl\prod_{i=1}^{n_1}\{x_i\}\time
Y_i\Br\times\Bl\prod_{j=1}^{m_1}Z_{\tau(j)}\Br\sub
\bl\sig_1\time X\br^{n_1}\time X^{m_1},\\
\Theta_2(\tau)&=\Bl\prod_{i=1}^{n_2}\{y_i\}\time
Y_{n_1+i}\Br\times\Bl\prod_{j=1}^{m_2}Z_{\tau(m_1+j)}\Br\sub
\bl\sig_2\time X\br^{n_2}\time X^{m_2},
\endalign
$$
that depend on the choice of $\tau\in S(m)$, a permutation of
$m$\ letters.
The subscheme we will work with is the intersection scheme
$$P_0\umo\Bl\{x_0,y_0\}\time \Del\Br\bigcap
\Bl \bl P^{n_1|m_1}\br\umo\bl\Theta_1(\tau)\br\time
\bl P^{n_2|m_2}\br\umo\bl\Theta_2(\tau)\br\Br,
\tag 3.11
$$
where
$$P_0: \cbnmoo\times\cbnmto\lra (\sigo\time\sigt)\times(X\time X)
$$
is defined as $P_0=(p_{0\sigo},q_{0\sigt})\times(p_{0X},q_{0X})$.
In the following, for $\ups_i\in \bl PGL(2)\time G\br^{n_i}\time G^{m_i}$,
we will denote by $\Theta_i(\tau)^{\ups_i}$\ the translation of
$\Theta_i(\tau)$\ by $\ups_i$.

\pro{Proposition 3.1}
Let the notation be as before and $\tau\in S(m)$. Then for general
$$(\ups_1,\ups_2)\in
\Bl \bl PGL(2)\time G\br^{n_1}\time G^{m_1}\Br\time
\Bl \bl PGL(2)\time G\br^{n_2}\time G^{m_2}\Br
\tag 3.12
$$
the scheme
$$P_0\umo\Bl\{x_0,y_0\}\time\Del\Br
\bigcap\Bl \bl P^{n_1|m_1}\br\umo\bl \Theta_1(\tau)^{\upsilon_1}\br
\times \bl P^{n_2|m_2}\br\umo\bl\Theta_2(\tau)^{\upsilon_2}\br\Br
\tag 3.13
$$
is discrete, Cohen-Macaulay and is contained in the fibers over
$W_{B_1}\ucirc\times W_{B_2}\ucirc$\ whose length is exactly
$$\align
\sum_{l=1}^p&\vphi_{B_1}\bl\zeta_l,\al_1\cdotss\al_{n_1}
\midd \be_{\tau(1)}\cdotss\be_{\tau(m_1)}\br\\
  &\quad \cdot\vphi_{B_2}\bl\tilde\zeta_l,\al_{n_1+1}
\cdotss\al_n \midd \be_{\tau(m_1+1)}\cdotss\be_{\tau(n)}\br
\endalign
$$
\endpro

\proof
We recall that by (3.8),
$$\dim \cbnmoo\times\cbnmto=
\codim\Theta_1(\tau)+\codim\Theta_2(\tau)
+\codim \{x_0,y_0\}\time\Del.
$$
Also, $\cbnmoo\times\cbnmto$\ is smooth over $W_{B_1}\ucirc\time
W_{B_2}\ucirc$\ and $\bl\sig_i\time X\br^{n_i}\time X^{m_i}$\
and $(\sigo\time\sigt)\time(X\time X)$\ are
homogeneous under $(PGL(2)\time G)^{n_i}\time G^{m_i}$\ and
$PGL(2)^2\time G^2$\ respectively. Hence for general
$\bl\upsilon_1,\ups_2,(\del_1,\del_2),(\gamma_1,\gamma_2)\br$\ in
$$\Bl \bl PGL(2)\time G\br^{n_1}\time G^{m_1}\Br\time
\Bl \bl PGL(2)\time G\br^{n_2}\time G^{m_2}\Br \time
PGL(2)^2\time G^2,
$$
the scheme
$$P_0\umo\Bl\{x_0^{\del_{1}},y_0^{\del_{2}}\}\time
\Del^{(\gamma_{1},\gamma_{2})}\Br
\bigcap\Bl \bl P^{n_1|m_1}\br\umo\bl \Theta_1(\tau)^{\upsilon_1}\br
\time \bl P^{n_2|m_2}\br\umo\bl\Theta_2(\tau)^{\upsilon_2}\br\Br
\tag 3.14
$$
is discrete, Cohen-Macaulay and is contained in the
fibers over $W_{B_1}\ucirc\times W_{B_2}\ucirc$, by Lemma 2.3.
To prove the first part of the proposition, we need to show that we can
achieve the same goal by only choosing $(\upsilon_1,\upsilon_2)$\
general while letting $\del_{1}$, $\del_{2}$, $\gamma_{1}$\
and $\gamma_{2}$\ be the identity elements of the respective groups.
Indeed, since $\sig\times X$\ is a $PGL(2)\time G$-variety, $C_B$\
admits a canonical $PGL(2)\time G$\ action that makes $(\ps,\px)\mh
C_B\to\sig\time X$\ equivariant. Hence
$$\bl p_{0\sigo},p_{0X}\br\times P^{n_1|m_1}:
C_{B_1}^{1+n_1+m_1}\lra
\bl\sig_1\time X\br\time\bl\sig_1\time X\br^{n_1}\time X^{m_1}
$$
is $PGL(2)\time G$\ equivariant with $PGL(2)$\ and $G$\
acting on both sides diagonally. Let $\mu_1$\ be the transformation
on $\bl\sig_1\time X\br^{n_1}\time X^{m_1}$\ that is induced by
$(\del_1,\gamma_1)$\
acting diagonally on it. We let $\mu_2$\ be the transformation
on $\bl\sig_2\time X\br^{n_2}\time X^{m_2}$\
defined similarly based on $(\del_2,\gamma_2)$. Then the scheme
$$P_0\umo\Bl\{x_0,y_0\}\time\Del\Br
\bigcap\Bl \bl P^{n_1|m_1}\br\umo
\bl \Theta_1(\tau)^{\mu_1\umo\cdot\upsilon_1}\br
\times \bl P^{n_2|m_2}\br\umo
\bl\Theta_2(\tau)^{\mu_2\umo\cdot\upsilon_2}\br\Br
$$
is canonically isomorphic to (3.14). Finally,
when $\bl\upsilon_1,\ups_2,(\del_1,\del_2),(\gamma_1,\gamma_2)\br$\
is a general element,
$$(\ups_1\pri,\ups_2\pri)=(\mu_1\umo\cdot\upsilon_1,
\mu_2\umo\cdot\upsilon_2)
$$
is a general element in (3.12) as well.
This proves the first part of the proposition.

Now we assume $(\ups_1,\ups_2)$\ is already a general element in (3.12).
Since the intersection scheme
(3.13) is discrete and Cohen-Macaulay, by [Fu, Proposition 7.1],
its length coincides with the degree of
$$\align
\Biggl(\Bl\bl P^{n_1|m_1}\br\sta\bl[\Theta_1(\tau)] & \dual\br
\times \bl P^{n_2|m_2}\br\sta
\bl \Theta_2(\tau)]\dual\br\Br\\
& \bigcup
\bl P_0\br\sta\bl[\{x_0,y_0\}]\dual\times[\Del]\dual\br\Biggr)
\Bigl[\cbnmoo\times\cbnmto\Bigr].
\tag 3.15
\endalign
$$
However,
$$\bl P_0\br\sta\bl[\{x_0,y_0\}]\dual\time[\Del]\dual\br
=\sum_{l=1}^p\bl p_{0\sigo}\sta e_0\times q_{0\sigt}\sta e_0\br\cup
\bl p_{0X}\sta\zeta_l\times q_{0X}\sta\tilde\zeta_l\br,
$$
because of the Kunneth decomposition (3.1).
Therefore (3.15) (as 0-cycle) is identical to
$$\align
\sum_{l=1}^p
\Bl\bl P^{n_1|m_1}& \br\sta\bl[\Theta_1(\tau)]\dual\br\bigcup
\bl p_{0\sigo}\sta e_0\cup p_{0X}\sta\zeta_l\br\Br\Bigl[\cbnmoo\Bigr]\\
& \times
\Bl\bl P^{n_1|m_1}\br\sta\bl[\Theta_2(\tau)]\dual\br\sta\bigcup
\bl q_{0\sigt}\sta e_0\cup q_{0X}\sta\zeta_l\br\Br
\Bigl[\cbnmto\Bigr].
\endalign
$$
Finally, we obviously have the identity
$$\align
\deg\, \Biggl(\Bl\bl P^{n_1|m_1}\br\sta\bl[\Theta_1(\tau)]\dual\br\bigcup&
\bl p_{0\sigo}\sta e_0\cup p_{0X}\sta\zeta_l\br\Br\Bigl[\cbnmoo\Bigr]
\biggr)
\\&=\vphi\bl
\zeta_l,\al_1\cdotss\al_{n_1}\midd\be_{\tau(1)}\cdotss\be_{\tau(m_1)}\br.
\endalign
$$
Thus, the degree of (3.15) is exactly the formula given
in the proposition as desired. This completes the proof of the
proposition.
\qed

In the following, we will denote the scheme in (3.11) by
$$\Int\{B_1,n_1,m_1,\tau\}
$$
that depends implicitly on the choice of
$$\zn\in\sig\ \text{and}\ \Yn,\Zm\sub X,
$$
which are assumed to be in general
position in the sense of the Corollary 2.7 and Proposition 3.1.
The main goal in the remainder of this section is to
construct an isomorphism
$$F:\tInt^{n_1:n_2}_B\{\zdot,\Ydot,\Zdot\}
\lra \coprod_{\Sb B_1+B_2=B,\, m_1+m_2=m\\ \tau\in S(m_1,m)\endSb}
\int_B\bl B_1,n_1,m_1,\tau\br,
$$
where $S(m_1,m)\sub S(m)$\ is the subset consisting of those
$\tau\in S(m)$\ such that
$$\tau(1)<\cdots<\tau(m_1)\ \text{and}\ \tau(m_1+1)<\cdots\tau(m).
$$
Note that the source and the target schemes of $F$\ are 0-schemes whose
lengths are given by the left and the right hand sides
of (3.2) respectively. Hence the composition law (3.2) follows if we
can construct and confirm the isomorphism of such an $F$.

We first explain how $F$\ is defined as a map between sets. Let
$$(S,\skl)\in\tInt^{n_1:n_2}_B\{\zdot,\Ydot,\Zdot\}
$$
be any point. By Corollary 2.7, $S$\ has exactly two irreducible
components. Thus the morphism
$\varphi_{Z_0}\mh S\to Z_0$\ must be an isomorphism. We identify $S$\
with $Z_0=\sigo\cup\sigt$\ by this isomorphism and let $f_i\mh
\sig_i\to X$\ be the morphism induced by inclusion $\sig_i\sub S
\sub Z_0\times X$\ and the projection $Z_0\times X\to X$.
Let $B_i=[f_i(\sig_i)]\in\aotilalg$. So $f_i\in\Mor(\sig_i,B_i)$.

We recall that
$$(h_1\!\circ p_Z)\umo(z_i)\cap Z_0\quad\text{and}\quad
(h_2\!\circ p_Z)\umo(z_i)\cap Z_0
$$
are closed points in $Z_0$\ away from its singular point.
Hence $s_i$\ lies in $\sigo$\ for $1\leq i\leq n_1$\ because it lies in
$(h_1\!\circ p_Z)\umo(z_i)$. For similar reason,
$s_{i}$\ lies in $\sigt$\ for $n_1+1\leq i\leq n$.
As to $s_i$\ with $i>n$, we know that each of them is contained
in the smooth locus of $S$\ and is contained
in one and only one $Z_0\time Z_j$\ among $1\leq j\leq m$, by (3) of Corollary
2.7.
Hence we can find unique $m_1\geq 0$\ and $\tau\in S(m_1,m)$\
such that
$$s_{n+\tau(j)}\in\sigo\ \text{for}\ 1\leq j\leq m_1\
\text{and}\ s_{n+\tau(m_1+j)}\in\sigt\ \ \text{for}\ 1\leq j\leq m_2.
\tag 3.16
$$
With these convention agreed, we assign
$$\align
x_i&=s_i,\ i=1\cdotss n_1\quad\text{and}\quad
x_{n_1+j}=s_{\tau(j)},\ j=1\cdotss m_1;\\
y_i&=s_{n_1+i},\ i=1\cdotss n_2\quad\text{and}\quad
y_{n_2+j}=s_{\tau(m_1+j)},\ j=1\cdotss m_2.
\tag 3.17
\endalign
$$
Then $f_1\mh\sigo\to X$\ is a morphism sending
$x_i$\ to $Y_i$\ for $i\leq n_1$\ and sending $x_{n_1+j}$\
to $Z_{\tau(j)}$\ for $j\leq m_1$. Similarly, $f_2\mh\sigt\to X$\
sends $y_i$\ to $Y_{n_1+i}$\ and $y_{n_2+j}$\ to $Z_{\tau(m_1+j)}$.
Finally, if we let $x_0\in\sigo$\ and $y_0\in\sigt$\ be as before
(the singular point of $Z_0$), then certainly
$$f_1(x_0)=f_2(y_0).
$$
Hence the graph $\Gamma_{f_1}\sub\sigo\times X$\ with tuple
$$\bigl\{(x_0,f_1(x_0))\cdotss(x_{n_1+m_1},f_1(x_{n_1+m_1})\bigr\}
$$
is a point in $\cbnmoo$\ and the graph $\Gamma_{f_2}$\ with
tuple $\{(y_i,f_2(y_i))\}$\ is a point in $\cbnmto$.
This pair of points assigns a point in
$\cbnmoo\times\cbnmto$\ that lies in
$$\Int\bl B_1,n_1,m_1,\tau\br,
$$
by straight forward inspection.
We denote this assignment by $F$, which for the moment is a
map between {\sl sets}\,:
$$F:\tInt^{n_1:n_2}_B\{\zdot,\Ydot,\Zdot\}\lra
\coprod\hugeindex \Int\bl B_1,n_1,m_1,\tau\br.
\tag 3.18
$$

\pro{Proposition 3.2}
Let the notation be as before. Suppose
$\zn\in Z_0$\ and $\Yn,\Zm\sub X$
are in general positions
as in Corollary 2.7 and Proposition 3.1, then
\roster
\item
the map $F$\ is one-to-one;
\item
The map $F$\ can be extended to a morphism $\tilF$\
of the respective schemes;
\item
$\tilF$\ is a local isomorphism at $w\in
\tInt^{n_1:n_2}_B\{\zdot,\Ydot,\Zdot\}$\ if $\pi\mh\h_B\to Z$\ is smooth
at $\bar w\in\h_B$, where $\bar w\in \h_B$\ lies under $w$.
\endroster
\endpro

\proof
(1) follows from the uniqueness of $m_1$\ and $\tau\in S(m_1,m)$\ which
follows from Corollary 2.7.
We now show that $F$\ can be extended to a morphism.
Since $\tInt^{n_1:n_2}_B\{\zdot,\Ydot,\Zdot\}$\
is zero-dimensional, it suffices to construct a morphism
$$\tilF_A: \spec A\lra \coprod\hugeindex\Int\bl B_1,n_1,m_1,\tau\br
$$
for any
$$\spec A\sub \tInt^{n_1:n_2}_B\{\zdot,\Ydot,\Zdot\}
\tag 3.19
$$
so that $\tilF_A$\ commutes with the base change,
where $A$\ is any local Artinian ring of residue field $K$.
Let $(S_A,s_{A,1}\cdotss s_{A,n+m})$\ be a flat family over $\spec A$\
induced by the inclusion (3.19), where
$$S_A\sub Z_0\times X\times\spec A
\tag 3.20
$$
and $s_{A,1}\cdotss s_{A,n+m}$\ are sections of $S_A\to\spec A$.
(By abuse of notation, we will also view $s_{A,i}$\ as a subscheme
of $S_A$\ that is isomorphic to $\spec A$\ under the projection
$S_A\to\spec A$.) (3.20) induces a morphism
$$\vphi_{Z_0}(A)\mh S_A\to Z_0\time\spec A
$$
after composing with projection to $Z_0\time\spec A$.
We first claim that $\vphi_{Z_0}(A)$\ is an
isomorphism. Indeed, let
$$\iota:\vphi_{Z_0}(A)\sta\OO_{Z_0\times\spec A}\lra \OO_{S_A}
\tag 3.21
$$
be the induced homomorphism. Let $S_0\sub S_A$\ be the scheme with
reduced scheme structure. By Corollary 2.7, $S_0$\ has exactly two
irreducible components, hence isomorphic to $Z_0$\ via $\varphi_{Z_0}\mh
S_0\to Z_0$\ (since $\xx(\OO_{S_0})=\xx(\OO_{Z_0})$).
On the other hand, since sheaves in (3.21) are flat over $A$,
(3.21) restricts to a homomorphism
$\vphi_{Z_0}\sta\OO_{Z_0}\to \OO_{S_0}$\ that is an
isomorphism since $Z_0\cong S_0$.
Hence $\iota$\ is injective because
both sheaves are flat over $A$. Further, because $\xx(\OO_{S_A})=
\xx(\OO_{Z_0\times\spec A})$, $\iota$\ must be an isomorphism. This
proves that $\vphi_{Z_0}(A)$\ is an isomorphism. Next after identifying
$S_A$\ with $Z_0\times\spec A$\ by $\vphi_{Z_0}(A)$, we obtain an
immersion
$$f_A: Z_0\time\spec A\sub Z_0\time X\time\spec A
$$
that provides us a pair of immersions
$$f_{A,i}: \sig_i\time\spec A\lra \sig_i\time X\time\spec A,\quad
i=1,2
$$
because $Z_0$\ is a union of $\sigo$\ and $\sigt$. Let
$$C_{A,i}\sub \sig_i\time X\time\spec A
$$
be the image of $f_{A,i}$. To proceed, we need to assign
sections of $C_{A,i}\to \spec A$\ based on $s_{A,i}$.
Let $m_1\geq 0$\ and $\tau\in S(m_1,m)$\ be the pair associated to
$$(s_A,s_{A,1}\cdotss s_{A,n+m})\otimes \spec K
$$
given at (3.16) and let
$$x_{A,i}, y_{A,j}:\spec A\sub S_A,\
1\leq i\leq n_1+m_1, 1\leq j\leq n_2+m_2
$$
be immersions derived by renaming $\{s_{A,i}\}_1^{n+m}$\
with rule specified in (3.17). Following this assignment, images of
$x_{A,i}$'s are contained in $C_{A,1}$\ and images of
$y_{A,j}$'s are contained in $C_{A,2}$.
Finally, we let $x_{A,0}\in C_{A,0}$\ be the image $f_{A,1}(x_0\time
\spec A)$\ and let $y_{A,0}\sub C_{A,2}$\ be the image
$f_{A,2}(y_0\time\spec A)$. Then
$$\bl C_{A,1},x_{A,0}\cdotss x_{A,n_1+m_1}\br
$$
defines a morphism
$$F_{A,1}: \spec A\lra \cbnmoo
$$
by the universality
of the Hilbert scheme $H_{B_1}\supseteq W_{B_1}$. Similarly, we obtain
a morphism
$$F_{A,2}: \spec A\lra \cbnmto
$$
by using $C_{A,2}$\ and $y_{A,j}$'s. Here $B_1$\ and $B_2$\ are
cycles in $A_1X$\ defined before based on two irreducible
components of $S_0\sub Z_0\time X$\ (before (3.16)).
Hence we obtain
$$F_A=F_{A,1}\times_A F_{A,2}:\spec A\lra \cbnmoo\times\cbnmto.
\tag 3.22
$$
Because $F_A$\ commutes with base change, which is apparent
from the construction, (3.22) defines a morphism
$$\tilF:\tInt^{n_1:n_2}_B\{\zdot,\Ydot,\Zdot\}\lra
\cbnmoo\times\cbnmto.
$$

We now show that the image lies in $\tInt^{n_1:n_2}_B\{\zdot,\Ydot,\Zdot\}$.
{}From the construction we see that (Im stands for image):
$$\cases
\varphi_X\bl\Im\{ x_{A,i}\}\br \sub Y_i\quad
\text{and}\quad \varphi_{\sigo}\bl\Im\{ x_{A,i}\}\br \sub \{z_i\},
& \text{for}\ 1\leq i\leq n_1;\\
\varphi_X\bl\Im\{ x_{A,n_1+j}\}\br \sub Z_{\tau(j)},
& \text{for}\ 1\leq j\leq m_1.
\endcases
$$
$$\cases
\varphi_X\bl\Im\{ y_{A,i}\}\br \sub Y_{n_1+i}\quad
\text{and}\quad \varphi_{\sigt}\bl\Im\{ y_{A,i}\}\br \sub \{z_{n_1+i}\},
& \text{for}\ 1\leq i\leq n_2;\\
\varphi_X\bl\Im\{ y_{A,n_2+j}\}\br \sub Z_{\tau(m_1+j)},
& \text{for}\ 1\leq j\leq m_2.
\endcases
$$
Also, it is straight forward to check that the image of the subscheme
$$x_{A,0}\times_A y_{A,0}\sub C_{A,1}\times_A\! C_{A,2}
$$
under the projection
$$C_{A,1}\times_A\! C_{A,2}
\sub \cbnmoo\times\cbnmto\mapright{P_0}
\bl \sigo\time\sigt\br\time\bl X\time X\br\lra X\time X
$$
is contained in $\Del$. Therefore, the morphism $\tilF$\ factor through
$$\coprod\hugeindex\Int\bl B_1,n_1,m_1,\tau\br
\sub\cbnmoo\times\cbnmto,
$$
which is an extension of $F$\ in (3.18). We
denote this extension by $\tilF$.

It remains to show that $\tilF$\ is a local isomorphism over where
$\pi$\ is smooth. Let
$$w\in\tInt^{n_1:n_2}_B\{\zdot,\Ydot,\Zdot\}
$$
be a smooth point of $\pi$\ and let $U$\ be the component of
$$\coprod\hugeindex\Int\bl B_1,n_1,m_1,\tau\br
\tag 3.23
$$
containing $\tilF(w)$. (Note (3.23) is a 0-dimensional scheme.)
Thus to show that $\tilF$\ is a local isomorphism near $w$\ it
suffices to construct a morphism
$$\widetilde{H}: U\lra \tInt^{n_1:n_2}_B\{\zdot,\Ydot,\Zdot\}
$$
such that $\widetilde{H}\circ \tilF$\ and
$\tilF\circ \widetilde{H}$\ are identity morphisms
on the components containing $w$\ and $\tilF(w)$\ respectively.
Let $U=\spec A$\ and let
$$\bl C_{A,1},x_{A,0}\cdotss x_{A,n_1+m_1}\br\times
\bl C_{A,2},y_{A,0}\cdotss y_{A,n_2+m_2}\br
\tag 3.24
$$
be the family over $\spec A$\ corresponding to
$$\spec A\sub
\coprod\hugeindex\Int\bl B_1,n_1,m_1,\tau\br.
$$
By Proposition 3.1, $C_{A,1}$\ and $C_{A,2}$\ are irreducible
and then
$$C_{A,1}\cong\sigo\time\spec A\ \text{and}\ C_{A,2}\cong\sigt\time\spec A
$$
via projections.
Let $S_A$\ be the scheme resulting from gluing $C_{A,1}$\
and $C_{A,2}$\ along $x_{A,0}\sub C_{A,1}$\ and
$y_{A,0}\sub C_{A,2}$.
Clearly, between {\it sets} we have a one-to-one map
$$f: S_A\lra Z_0\times X\times\spec A
\tag 3.25
$$
that is induced by the closed immersion
$$f_1: C_{A,1}\lra \sigo\time X\time\spec A,\quad
f_2: C_{A,2}\lra \sigt\time X\time\spec A.
$$
To show that (3.25) actually comes from a closed immersion (flat over
$\spec A$), we first need to define a homomorphism
$$f\sta\OO_{Z_0\times X\times\spec A}\lra \OO_{S_A}.
\tag 3.26
$$
Recall that $\OO_{S_A}$\ is defined by the exact sequence
$$0\lra\OO_{S_A}\lra \OO_{C_{A,1}}\oplus \OO_{C_{A,2}}\lra A\lra 0
$$
as $\OO_{S_A}$-modules (by viewing $C_{A,1}, C_{A,2}\sub S_A$\ and
identifying $A=\OO_{x_{A,0}}$\ and $A=\OO_{y_{A,0}}$)
and $\OO_{Z_0\times X\times\spec A}$\ belongs to the exact sequence
$$0\lra\OO_{Z_0\times X\times\spec A}\lra \OO_{\sigo\times X\times\spec A}
\oplus \OO_{\sigt\times X\times\spec A}\lra \OO_{X\times \spec A}\lra 0.
$$
Hence to find (3.26) we only need to show that the composition
$$f\sta \OO_{Z_0\times X\times\spec A}\lra
f_1\sta \OO_{\sigo\times X\times\spec A}
\oplus f_2\sta \OO_{\sigt\times X\times\spec A}\lra
\OO_{C_{A,1}}\oplus \OO_{C_{A,2}}\lra A
$$
is trivial, which is true because of the condition
$$\Int\bl B_1,n_1,m_1,\tau\br\sub
P_0\umo(\{x_0,y_0\}\times\Del).
$$
Therefore, we obtain the homomorphism (3.26).
This homomorphism is surjective because
$$f_1\sta\OO_{\sigo\times X\times\spec A}\lra \OO_{C_{A,1}}\
\text{and}\ f_2\sta\OO_{\sigt\times X\times\spec A}\lra \OO_{C_{A,2}}
$$
are surjective and henceforth induce a morphism
$$g:\spec A\lra \h_{B,0}.
\tag 3.27
$$
We claim that $g$\ factor through $\w_{B,0}$. Indeed, by our
construction of $\tilF$\ and $g$, $\bar w=g(\spec K)\in \h_{B,0}$\ is
under $w\in\c_{B,0}^{n+m}$\ via the projection $\c_{B,0}^{n+m}
\to \w_{B,0}$. Hence $\bar w$\ is contained in $\w_{B,0}$.
However, $\bar w$\ is a smooth point of $\pi\mh\h_B\to V$\ by assumption.
Therefore $g$\ must factor through
$\w_{B,0}\sub\h_{B,0}$\ because a
neighborhood of $g(\spec K)\in\h_{B,0}$\ is contained in $\w_{B,0}$.

Now we construct $\widetilde{H}$. Let $\tau\in S(m_1,m)$\ be the permutation
associated to $w$\ specified in (3.16) and let
$$s_{A,i}:\spec A\lra S_A,\ i=1\cdotss n+m
$$
be sections derived from $\{x_{A,i}\}^{n_1+m_1}$\ and
$\{y_{A,i}\}^{n_2+m_2}$\ by the rule (3.17). Since $g$\ factor through
$\w_{B,0}\sub\h_{B,0}$, the tuple
$(S_A,s_{A,1}\cdotss s_{A,n+m})$\ defines a morphism
$$\widetilde{H}:\spec A\lra \c_{B,0}^{n+m}
$$
that factor through $\bl P^{n_1:n_2|m}\br\umo\bl\widetilde Y\br$,
by our choice of $S_{A,i}$'s. This gives $\widetilde{H}$.

{}From the construction of $\tilF$\ and $\widetilde{H}$, it is straight forward
to check that $\tilF\circ \widetilde{H}$\ and
$\widetilde{H}\circ \tilF$\ are identities on the respective components.
This completes the proof of the proposition.
\qed

With all these prepared, the proof of the our main theorem
is a matter of formality.

\noindent
{\sl Proof of Theorem 0.2}.
We need to show that for any $B\in\aotilalg$\ and
$n=n_1+n_2$\ with $n_1,n_2\geq 2$,
$$\align
&\phinm(\aln\midd\bem)\\
=&\sum_{\Sb B=B_1+B_2,\,\tau\in S(m)\\
1\leq k\leq m,\,1\leq l\leq p\endSb}
    {1\over k!(m-k)!}\cdot
 \vphi_{B_1}\Bigl(\zeta_l,\al_1\cdotss \al_{n_1}\midd
      \be_{\tau(1)}\cdotss\be_{\tau(k)}\Bigr)\\
    &\qquad\qquad\qquad\times
\vphi_{B_2}\left(\tilde\zeta_l,\al_{n_1+1}\cdotss \al_n
\midd \be_{\tau(k+1)}\cdotss \be_{\tau(m)}\right).
\tag 3.28
\endalign
$$

We first remark that since $\vphi_B$\ is a group homomorphism,
it suffices to consider the classes $\aln,\bem\in A\sta X$\ that
are Poincare dual to subvarieties $\Yn$, $\Zm\sub X$.
Also, by dimension reason, (3.28) is non-trivial only when
$$\sum_{i=1}^n\codim Y_i +\sum_{j=1}^m\codim Z_j=
\rho(B)+m,
$$
which we will assume in the remainder of the proof.

By using the degeneration $\w_B$\ and the corresponding family
$\c_B$, we can express the value
$$\vphi_B(\aln\midd\bem)
$$
as the degree of the zero cycle (see (3.7))
$$\tilde\vpsi_B^{n_1:n_2}\bigl(\aln\midd\bem\bigr)
\bigl[\c_{B,0}^{n+m}\bigr].
$$
We choose $\zn\in\sig$\ and $\Yn,\Zm\sub X$\ in general
position (among their translations by the respective groups) as
specified in Corollary 2.7 and Proposition 3.1 and let
$\widetilde Y$\ be the subvariety defined before (before (3.9)).
Following Proposition 3.1 and 3.2, we know that
points in
$$\bl P^{n_1:n_2|m}\br\umo\bl\widetilde Y\br\cap
\c_{B,0}^{n+m}
\tag 3.29
$$
are represented by
$$(S,s_1\cdotss s_{n+m})
$$
such that $S\sub Z_0\time X$\ have exactly two irreducible components
and all $s_i\in S$\ are contained in the smooth locus of $S$.
Hence $\c_B^{n+m}$\ is smooth near (3.29), by Lemma 2.5
and [Ko, \S1.2]. Therefore, by Lemma 2.3,
$$\bl P^{n_1:n_2|m}\br\umo\bl\widetilde Y\br
\tag 3.30
$$
is one-dimensional, Cohen-Macaulay (smooth
when $\char K=0$) near $\c_{B,0}^{n+m}$, after making $z_i\in\sig$\
and $Y_i,Z_j\sub X$\ in general position if necessary.
Finally, we conclude that
(3.30) is flat over $V$\ near $0\in V$\ because the map
$$\tilF:\bl P^{n_1:n_2|m}\br\umo\bl\widetilde Y\br\cap
\c_{B,0}^{n+m}\lra
\coprod\hugeindex\Int\bl B_1,n_1,m_1,\tau\br
$$
is one-one and the target is 0-dimensional (see Proposition 3.2).
This shows that the intersection scheme
$$\tInt^{n_1:n_2}_B\{\zdot,\Ydot,\Zdot\}
=\bl P^{n_1:n_2|m}\br\umo\bl\widetilde Y\br\cap
\c_{B,0}^{n+m}
$$
is 0-dimension and Cohen-Macaulay and by [Fu, Proposition 7.1 and
Corollary 17.4] and (3.6),
$$\align
\vphi_B(\aln\midd\bem)
=&\deg \Bl\tilde\vpsi_B^{n_1:n_2}(\aln\midd\bem)[\c_{B,0}^{n+m}]\Br\\
=&\text{length}\Bl\tInt^{n_1:n_2}_B\{\zdot,\Ydot,\Zdot\}\Br.
\tag 3.31
\endalign
$$

Next, we turn our attention to the morphism
$$\tilF:\tInt^{n_1:n_2}_B\{\zdot,\Ydot,\Zdot\}\lra
\coprod\hugeindex \Int\bl B_1,n_1,m_1,\tau\br.
$$
We first show that $\tilF$\ is an isomorphism. By Proposition 3.2
it suffices to show that $\tilF$\ is surjective and the condition in (3)
of Proposition 3.2 holds everywhere. Let
$$\bigl\{(C_1,x_0\cdotss x_{n_1+m_1}), (C_2,y_0\cdotss y_{n_2+m_2})
\bigr\}
$$
be any point in $\Int\bl B_1,n_1,m_1,\tau\br$. By Corollary 2.7, both
$C_1$\ and $C_2$\ are irreducible. As we did in the proof of
Proposition 3.2, we get a curve $S\sub Z_0\time X$\ that is the union of
the image of $C_1$\ and $C_2$\ in $Z_0\time X$.
We also obtain ordered points $s_1\cdotss s_{n+m}\in S$\ that
is the result of renaming $x_1\cdotss x_{n_1+m_1}$\ and $y_1\cdotss
y_{n_2+m_2}$\ according to the rule (3.17) that depends on $\tau\in
S(m_1,m)$. From the construction, it is clear that the tuple
$(S,s_1\cdotss s_{n+m})$\ represents a point in
$\tInt^{n_1:n_2}_B\{\zdot,\Ydot,\Zdot\}$\ when the point $w\in \h_{B,0}$\
associated to $S$\ is in $\w_{B,0}$. The reason $w$\ does belong
to $\w_{B,0}$\ is as follows:
By (1) of Lemma 2.5, $\h_B$\ is smooth at $w$. Also by (2) of
Lemma 2.5, $\pi\mh \h_B\to V$\ is smooth at $w$. Let $U$\ be the
component of $\h_B$\ containing $w$. Then general points of $U$\
must be irreducible curves in $\sig\time X$, since $S$\ has exactly
two irreducible components. Therefore, $U$\
is contained in the closure of $W_B\ucirc\time(V-0)$. This proves that
$w\in\w_{B,0}$\ and consequently, $(S,s_1\cdotss s_{n+m})$\
represents a point $\bar w\in \tInt^{n_1:n_2}_B\{\zdot,\Ydot,\Zdot\}$.
By our construction, we certainly have $\tilF(\bar w)=w$. Hence
$\tilF$\ is surjective. This argument also shows that $\pi\mh\h_B\to V$\
is smooth at all $\bar w\in \tInt^{n_1:n_2}_B\{\zdot,\Ydot,\Zdot\}$\.
Therefore by Proposition 3.2, $\tilF$\ is an isomorphism.
Combined with (3.31), we have
$$\align
&\phinm(\aln\midd\bem)\\
=&\sum_{\Sb B=B_1+B_2\\ 0\leq k\leq m,\,\tau\in S(k,m)\endSb}
\text{length}\, \Bl\Int\bl B_1,n_1,m_1,\tau\br\Br\\
=&\sum_{\Sb B=B_1+B_2,\,0\leq k\leq m\\1\leq l\leq p,\,\tau\in S(k,m)
\endSb}
 \vphi_{B_1}\Bigl(\zeta_l,\al_1\cdotss \al_{n_1}\midd
      \be_{\tau(1)}\cdotss\be_{\tau(k)}\Bigr)\\
    &\qquad\qquad\qquad\qquad\qquad\cdot
\vphi_{B_2}\left(\tilde\zeta_l,\al_{n_1+1}\cdots \al_n
\midd \be_{\tau(k+1)}\cdotss \be_{\tau(m)}\right)\\
=&\sum_{\Sb B=B_1+B_2,\,0\leq k\leq m\\1\leq l\leq p,\,\tau\in S(m)
\endSb}
{1\over k!\,(m-k)!}\vphi_{B_1}\Bigl(\zeta_l,\al_1\cdotss \al_{n_1}\midd
      \be_{\tau(1)}\cdotss\be_{\tau(k)}\Bigr)\\
    &\qquad\qquad\qquad\qquad\qquad\cdot
\vphi_{B_2}\left(\tilde\zeta_l,\al_{n_1+1}\cdots \al_n
\midd \be_{\tau(k+1)}\cdotss \be_{\tau(m)}\right).
\endalign
$$
Here the second identity follows from Proposition 3.1 and the last
identity follows from
the symmetry of $\vphi_B(\aln\midd\bem)$\ among $\al_i$'s
and among $\be_j$'s.
This completes the proof of the composition law.
\qed

\vskip30pt
\noindent
Mathematics Department,
University of California, Los Angeles, CA 90024

\noindent
jli\@math.ucla.edu;

\vskip15pt
\noindent
Courant Institute of Mathematics,
New York University, New York NY 10012

\noindent
tiang\@math1.cims.nyu.edu.

\bye